\documentclass[aop,seceqn,citesort,MSNbibl,dvips]{arximspdf}

\doi{10.1214/09-AOP496}
\volume{38}
\issue{2}
\pubyear{2010}
\firstpage{794}
\lastpage{840}

\makeatletter

\newproclaim{Remark}{Remark}[section]
\newtheorem{Theorem}{Theorem}[section]
\newtheorem{Lemma}{Lemma}[section]
\newproclaim{Definition}{Definition}[section]
\newtheorem{Proposition}{Proposition}[section]

\makeatother

\begin{document}
\begin{frontmatter}

\title{Backward SDEs with constrained jumps and quasi-variational inequalities}
\runtitle{BSDEs with constrained jumps and QVIs}

\begin{aug}
\author[A]{\fnms{Idris} \snm{Kharroubi}\ead[label=e1]{kharroubi@ensae.fr}},
\author[B]{\fnms{Jin} \snm{Ma}\thanksref{t1}\ead[label=e2]{jinma@usc.edu}},
\author[A]{\fnms{Huy\^{e}n} \snm{Pham}\corref{}\ead[label=e3]{pham@math.jussieu.fr}} and
\author[B]{\fnms{Jianfeng} \snm{Zhang}\thanksref{t2}\ead[label=e4]{jianfenz@usc.edu}}
\runauthor{Kharroubi, Ma, Pham and Zhang}
\affiliation{PMA Universit{\'{e}} Paris 7 and CREST, University of
Southern California, PMA~Universit{\'{e}} Paris 7, CREST and IUF
and University of Southern California}
\address[A]{I. Kharroubi\\
H. Pham\\
Laboratoire de Probabilit\'{e}s\\
\quad et Mod\`{e}les Al\'{e}atoires\\
CNRS, UMR 7599\\
Universit{\'{e}} Paris 7\\
Case courrier 7012\\
2, Place Jussieu\\
75251 PARIS Cedex 05\\
France \\
\printead{e1}\\
\phantom{E-mail: }\printead*{e3}} 
\address[B]{J. Ma\\
J. Zhang\\
Department of Mathematics\\
University of Southern California\\
3620 Vermont Avenue, KAP 108\\
Los Angeles, California 90089\\
USA\\
\printead{e2}\\
\phantom{E-mail: }\printead*{e4}}
\end{aug}

\thankstext{t1}{Supported in part by NSF Grant 0505427.}
\thankstext{t2}{Supported in part by NSF Grant 0631366.}

\pdfauthor{Idris Kharroubi, Jin Ma, Huyen Pham, Jianfeng Zhang}

\received{\smonth{6} \syear{2008}}
\revised{\smonth{3} \syear{2009}}

%
\begin{abstract}
We consider a class of backward stochastic differential equations
(\mbox{BSDEs}) driven by Brownian motion and Poisson random measure, and
subject to constraints on the jump component. We prove the existence
and uniqueness of the minimal solution for the BSDEs by using a
penalization approach. Moreover, we show that under mild conditions
the minimal solutions to these constrained BSDEs can be
characterized as the unique viscosity solution of quasi-variational
inequalities (QVIs), which leads to a probabilistic representation
for solutions to QVIs. Such a representation in particular gives a
new stochastic formula for value functions of a class of impulse
control problems. As a direct consequence, this suggests a numerical
scheme for the solution of such QVIs via the simulation of the
penalized BSDEs.
\end{abstract}

%

\end{frontmatter}

\section{Introduction and summary}\label{sec1}

Consider a parabolic quasi-variational inequality (QVI for short)
of the following form:
%
%
\begin{equation}\label{IQV}
\cases{\min\biggl[ - \dfrac{\partial v}{\partial t} -
\mathcal{L}v - f , v - \mathcal{H}v \biggr] = 0, &\quad
on $[0,T)\times\mathbb{R}^d$,\cr
v(T,\cdot) = g, &\quad on $\mathbb{R}^d$,}
\end{equation}
where $\mathcal{L}$ is the second-order local operator
%
%
\begin{equation}\label{localL}
\mathcal{L}v(t,x) = \langle b(x),D_x v(t,x)\rangle+ \tfrac
{1}{2}\operatorname{tr}(\sigma
\sigma^{\intercal}(x)D_x^2v(t,x))
\end{equation}
and $\mathcal{H}$ is the nonlocal operator
%
%
\begin{equation}\label{nonlocalH} \mathcal{H}v
(t,x) = \sup_{e\in E} \bigl[ v\bigl(t,x+\gamma(x,e)\bigr) + c(x,e) \bigr].
\end{equation}
In the above, $D_x v$ and $D_x^2v$ are the partial gradient and
the Hessian matrix of $v$ with respect to its second variable $x$,
respectively; $^{\intercal}$ stands for the transpose;
$\langle\cdot,\cdot\rangle$ denotes the scalar product in
$\mathbb{R}^d$; $\mathbb{S}^d$ is
the set of all symmetric $d\times d$ matrices; and $E$ is some
compact subset of $\mathbb{R}^q$.

It is well known (see, e.g.,
\cite{benlio84}) that the QVI (\ref{IQV}) is the dynamic programming
equation associated to the impulse control problems whose value
function is defined by
%
%
\begin{eqnarray}\label{impulse}
v(t,x) &=&
\sup_{\alpha=(\tau_i,\xi_i)_i}\mathbf{E} \biggl[
g(X_T^{t,x,\alpha}) +
\int_t^T f(X_s^{t,x,\alpha}) \,ds + \sum
_{t < \tau_i\leq T}
c(X_{\tau_i^-}^{t,x,\alpha},\xi_i) \biggr].\hspace*{-33pt}
\end{eqnarray}
More precisely, given a filtered probability space $(\Omega,\mathcal
{F},\mathbf{P},
\mathbb{F})$ where
$\mathbb{F}=\{\mathcal{F}_t\}_t$, we define an impulse control
$\alpha$ as a double
sequence $(\tau_i,\xi_i)_i$ in which $\{\tau_i\}$ is an increasing
sequence of $\mathbb{F}$-stopping times, and each $\xi_i$ is an
$\mathcal{F}_{\tau_i}$-measurable random variable taking values in
$E$. For
each impulse control $\alpha=(\tau_i,\xi_i)_i$, the controlled
dynamics starting from $x$ at time $t$, denoted by
$X^{t,x,\alpha}$, is a c\`{a}dl\`{a}g process satisfying
the following SDE:
%
%
\begin{equation}\label{controlX}
X_s^{t,x,\alpha} = x + \int_t^s
b(X_u^{t,x,\alpha}) \,du + \int_t^s \sigma
(X_u^{t,x,\alpha})\,
dW_u + \sum_{t<\tau_i\leq s} \gamma
(X_{\tau_{i}^-}^{t,x,\alpha},\xi_{i}),\hspace*{-33pt}
\end{equation}
where $W$ is a $d$-dimensional $\mathbb{F}$-Brownian motion. In other
words, the controlled process $X^{t,x,\alpha}$ evolves according to
a diffusion process between two successive intervention times
$\tau_i$ and $\tau_{i+1}$, and at each decided intervention time
$\tau_i$, the process jumps with size $\Delta
X_{\tau_i}^{t,x,\alpha}:=
X_{\tau_i}^{t,x,\alpha}-X_{\tau_i^-}^{t,x,\alpha}=\gamma(X_{\tau
_{i}^-}^{t,x,
\alpha},\xi_i)$.

We note that the impulse control problem (\ref{impulse}) may be
viewed as a sequence of optimal stopping problems combined with
jumps in state due to impulse values. Moreover, the QVI (\ref{IQV})
is the infinitesimal derivation of the dynamic programming
principle, which means that at each time, the controller may decide
either to do nothing and let the state process diffuse, or to make
an intervention on the system via some impulse value. The former is
characterized by the linear PDE in~(\ref{IQV}), while the latter is
expressed by the obstacle (or reflected) part in (\ref{IQV}). From
the theoretical and numerical point of view, the main difficulty of
the QVI (\ref{IQV}) lies in that the obstacle contains the solution
itself, and it is nonlocal [see (\ref{nonlocalH})] due to the jumps
induced by the impulse control. These features make the classical
approach of numerically solving such impulse control problems
particular challenging.
%

An alternative method to attack the QVI (\ref{IQV}) is to
find the probabilistic representation of the solution using the
backward stochastic differential equations (BSDEs), namely the
so-called nonlinear Feynman--Kac formula. One can then hope to use
such a representation to derive a direct
numerical procedure for the solution of QVIs, whence the impulse
control problems. The idea is
the following. We consider a Poisson random measure $\mu(dt,de)$
on $\mathbb{R}_+\times E$ associated to a marked point process
$(T_i,\zeta_i)_i$. Assume that $\mu$ is independent of $W$ and has
intensity $\lambda(de)\,dt$, where $\lambda$ is a finite measure on
$E$.
Consider a (uncontrolled) jump-diffusion process
%
%
\begin{equation}
\label{jumpdiff}
X_s = X_0 + \int_0^s b(X_u) \,du + \int_0^s
\sigma(X_u) \,dW_u + \sum_{T_i\leq s} \gamma(X_{T_i^-},\zeta_i).
\end{equation}
Assume that $v$ is a ``smooth'' solution to (\ref{IQV}), and define
$Y_t = v(t,X_t)$. Then, by It\^{o}'s formula, we have
%
%
\begin{eqnarray}\label{BSDEgen0pre}
Y_t &=& g(X_T) + \int_t^T f(X_s) \,ds + K_T - K_t
- \int_t^T \langle Z_s, dW_s\rangle\nonumber\\[-8pt]\\[-8pt]
&&{} - \int_t^T\int_E \bigl( U_s(e) - c(X_{s^-},e) \bigr) \mu(ds,de),\nonumber
\end{eqnarray}
where\vspace*{1pt} $Z_t=\sigma^{\intercal}(X_{t^-})D_xv(t,X_{t^-})$, $ U_t(e)=
v(t,X_{t^-}+\gamma(X_{t^-},e))- v(t,X_{t^-}) + c(X_{t^-},e)$
and $K_t = \int_0^t (-\frac{\partial v}{\partial t}-\mathcal{L}v -
f)(s,X_s) \,ds$. Since $v$
satisfies (\ref{IQV}), we see that $K$ is a continuous (hence,
predictable), nondecreasing process and $U$ satisfies the
constraint
%
%
\begin{equation}
\label{conspre}
- U_t(e) \geq0.
\end{equation}
The idea is then to view (\ref{BSDEgen0pre}) and (\ref{conspre}) as a
BSDE
with jump constraints, and we expect to retrieve $v(t,X_t)$ by
solving the ``minimal'' solution $(Y,Z,U,K)$ to this constrained
BSDE.

We can also look at the BSDE above slightly differently.
Let us denote $d{\bar K}_t = dK_t - \int_E U_s(e) \mu(dt,de)$, $t\ge
0$. Then $\bar K$ is still a nondecreasing process, and
equation (\ref{BSDEgen0pre}) can now be rewritten as
%
%
\begin{eqnarray}\label{Yfor}
Y_t & = & g(X_T) + \int_t^T f(X_s) \,ds + \int_t^T\int_E c(X_{s^-},e)
\mu(ds,de)
\nonumber\\[-8pt]\\[-8pt]
&&{} - \int_t^T \langle Z_s, dW_s\rangle+ \bar K_T - \bar K_t.\nonumber
\end{eqnarray}
We shall prove that
$v(t,X_t)$ can also be retrieved by looking at the minimal solution
$(Y,Z,{\bar K})$ to this BSDE. In fact, the following relation holds
(assuming $t=0$):
%
%
\begin{eqnarray}\label{super}\quad
v(0,X_0) & = & \inf\biggl\{ y \in\mathbb{R}\dvtx\exists
Z, y + \int_0^T \langle Z_s, dW_s\rangle\nonumber\\
&&\hspace*{138pt}\qquad\hspace*{-94.5pt}{}\geq g(X_T) + \int_0^T
f(X_s) \,ds \\
&&\hspace*{65.4pt}{} + \int_0^T \int_{E}c(X_{s^-},e) \mu(ds,de) \biggr\}.
\nonumber
\end{eqnarray}
Notice that (\ref{super}) also has a financial interpretation. That is,
$v(0,x)$ is the minimal capital allowing to
superhedge the payoff $\Pi_T(X)=g(X_T)+ \int_0^T f(X_s) \, ds +
\int_0^T c(X_{s^-},e)\mu(ds,de)$ by trading only the asset $W$.
Here, the market is obviously incomplete, since
the jump part of the underlying asset $X$ is not hedgeable.
This connection between the impulse control problem (\ref{impulse}) and
the stochastic target problem defined by the r.h.s. of (\ref{super}) was
originally proved in Bouchard \cite{bou06}.

Inspired by the above discussion, we now introduce the following
general BSDE:
%
%
\begin{eqnarray}\label{BSDEgen0}\hspace*{28pt}
Y_t & = & g(X_T) + \int_t^T f(X_s,Y_s,Z_s)\, ds + K_T - K_t
- \int_t^T \langle Z_s, dW_s \rangle\nonumber\\[-8pt]\\[-8pt]
& &{} - \int_t^T\int_E \bigl( U_s(e)
- c(X_{s^-},Y_{s^-},Z_{s},e) \bigr) \mu(ds,de),\qquad 0 \leq
t \leq T, \nonumber
\end{eqnarray}
with constraints on the jump component in
the form
%
%
\begin{equation}
\label{hcons0}
h(U_t(e)) \geq0,\qquad
\forall e \in E, 0 \leq t \leq T,
\end{equation}
where $h$ is a given nonincreasing function. The solution to the
BSDE is a quadruple $(Y,Z,U,K)$ where, besides the usual
component $(Y, Z, U)$,
the fourth component $K$ is a nondecreasing, c\`{a}dl\`{a}g, adapted
process, null at zero, which makes the constraint (\ref{hcons0})
possible.
We note that without the constraint (\ref{hcons0}), the BSDE with
$K=0$
was studied by Tang and Li \cite{tanli94} and Barles, Buckdahn and
Pardoux \cite{barbucpar97}. However, with the presence of the
constraint, we may not have the uniqueness of the solution.
We thus look only for the minimal solution $(Y,Z,U,K)$, in the sense
that
for any\vspace*{1pt} other solution $(\tilde Y,\tilde Z,\tilde U,\tilde K)$
satisfying (\ref{BSDEgen0}) and (\ref{hcons0}), it must hold that $Y
\leq\tilde Y$. Clearly, this BSDE is a generalized version of
(\ref{BSDEgen0pre}) and (\ref{conspre}), where the functions $f$ and $c$
are independent of $y$ and $z$, and
$h(u)=-u$.

We can also consider the counterpart of (\ref{Yfor}), namely finding
the minimal solution $(Y,Z,K)$ of the BSDE
%
%
\begin{eqnarray}\label{BSDE1}
Y_t &=& g(X_T) + \int_t^T f(X_s,Y_s,Z_s) \, ds\nonumber\\
&&{}  + \int_t^T\int_E c(X_{s^-},Y_{s^-},Z_{s},e)
\mu(ds,de) \\
&&{} - \int_t^T \langle Z_s, dW_s\rangle+ K_T - K_t,\qquad
0 \leq t \leq T.\nonumber
\end{eqnarray}
It is then conceivable, as we shall prove, that this problem is a
special case of (\ref{BSDEgen0}) and (\ref{hcons0}) with $h(u)=-u$.

It is worth noting that if the generator $f$ and the cost function
$c$ do not depend on $y,z$, which we refer to as the impulse
control case, the existence of a minimal solution to the constrained
BSDEs (\ref{BSDEgen0pre}) and (\ref{conspre}) may be directly obtained by
supermartingale decomposition method in the spirit of El Karoui and
Quenez \cite{elkque95} for the dual representation of the
super-replication cost of $\Pi_T(X)$. In fact, the results could be
extended easily to the case where $f$ is linear in $z$, via a simple
application of the Girsanov transformation. In our general case,
however, we
shall follow a penalization method, as was done in El Karoui et al.
\cite{elketal97}. Namely, we construct a suitable sequence
$(Y^n,Z^n,U^n,K^n)$ of BSDEs with jumps, and prove that it converges
to the minimal solution that we are looking for. This is achieved as
follows. We first show the convergence of the sequence
$(Y^n)$ by relying on comparison results for BSDEs with jumps, see
\cite{roy06}. The proof of convergence of the components
$(Z^n,U^n,K^n)$ is more delicate, and is obtained by using a weak
compactness argument due to Peng \cite{pen99}.

Our next task of this paper is to relate the minimal solution to the
BSDE with constrained jumps to the viscosity solutions to the
following general QVI:
%
%
\begin{equation}\label{QVIgen0}
\min\biggl[ - \frac{\partial v}{\partial t} - \mathcal{L}v - f(\cdot
,v,\sigma^{\intercal}D_x v) ,
h( \mathcal{H}v -
v) \biggr] = 0,
\end{equation}
where $\mathcal{H}$ is the nonlocal semilinear operator
\[
\mathcal{H}v (t,x)=
\sup_{e\in E} \bigl[ v\bigl(t,x+\gamma(x,e)\bigr) + c(x,v(t,x),\sigma
^{\intercal}(x)D_x
v(t,x),e) \bigr].
\]
Under suitable assumptions, we shall also prove the uniqueness of
the viscosity solution,
leading to a new probabilistic representation for this parabolic
QVI.

We should point out that BSDEs with constraints have been studied by
many authors. For example,
El Karoui et al. \cite{elketal97} studied the reflected BSDEs, in
which the component $Y$ is forced to stay above a given obstacle;
Cvitanic, Karatzas and Soner \cite{cvikarson98}, and Buckdahn
and Hu \cite{buchu98} considered the case where the constraints are
imposed on the component $Z$. Recently, Peng \cite{pen99} (see also
\cite{penxu07}) studied the the general case where constraints are
given on both $Y$ and $Z$, which relates these constrained BSDEs to
variational inequalities. The main feature of this work is to
consider constraints on the jump component ($U$) of the solution,
and to relate these jump-constrained BSDEs to quasi-variational
inequalities. On the other hand, the classical approach in the
theory and numerical approximation of impulse control problems and
QVIs is to consider them as obstacle problems and iterated optimal
stopping problems. However, our penalization procedure for
jump-constrained BSDEs suggests a noniterative approximation scheme
for QVIs, based on the simulation of the BSDEs, which, to our best
knowledge, is new.

The rest of the paper is organized as follows: in Section \ref{sec2}, we
give a detailed formulation of BSDEs with constrained jumps, and
show how it includes problem (\ref{BSDE1}) as special case.
Moreover, in the special case of impulse control, we directly
construct and show the existence of a minimal solution. In
Section \ref{sec3}, we develop the penalization approach for studying the
existence of a minimal solution to our constrained BSDE for
general $f$, $c$ and $h$. We show in Section \ref{secrelQVI} that the minimal
solution to this constrained BSDE provides a probabilistic
representation for the unique viscosity solution to a parabolic
QVI.
Finally, in Section \ref{secsuff}, we
provide some examples of sufficient conditions under which our
general assumptions are satisfied.

\section{BSDEs with constrained jumps}\label{sec2}

\subsection{General formulation} \label{secgeneral}

Throughout this paper, we assume that $(\Omega,\mathcal{F},\mathbf
{P})$ is a
complete probability space on which are defined a $d$-dimensional
standard Brownian motion $W = (W_t)_{t\geq0}$, and a Poisson
random measure $\mu$ on $\mathbb{R}_+\times E$, where $E$ is a
compact set
of $\mathbb{R}^q$, endowed with its Borel field $\mathcal{E}$. We
assume that the
Poisson random measure $\mu$ is independent of $W$, and has the
intensity measure $\lambda(de)\,dt$ for some finite measure $\lambda$
on $(E, \mathcal{E})$. 
We set $\tilde\mu(dt,de) = \mu(dt,de)-\lambda(de)\,dt$, the
compensated measure associated to $\mu$; and denote by $\mathbb{F} =
(\mathcal{F}_t)_{t\geq0}$ the augmentation of the natural filtration
generated by $W$ and $\mu$, and by $\mathcal{P}$ the $\sigma
$-algebra of
predictable subsets of $\Omega\times[0,T]$.

Given Lipschitz functions $b \dvtx\mathbb{R}^d \rightarrow
\mathbb{R}^d$,
$\sigma\dvtx\mathbb{R}^d \rightarrow\mathbb{R}^{d\times d}$,
and a measurable
map $\gamma\dvtx\mathbb{R}^d\times E \rightarrow\mathbb
{R}^d$, satisfying for
some positive constants $C$ and $k_\gamma$,
\[
{\sup_{e\in E}} |\gamma(x,e)| \leq C\quad \mbox{and}\quad
{\sup_{e\in E}}|
\gamma(x,e)- \gamma(x',e)| \leq k_\gamma|x-x'|,\qquad
x,x'\in\mathbb{R}^d,
\]
we consider the forward SDE:
%
%
\begin{equation}\label{defXjump}
dX_s = b(X_s)\,
ds + \sigma(X_s) \,dW_s + \int_E \gamma(X_{s^-},e) \mu(ds,de).
\end{equation}
Existence and uniqueness of (\ref{defXjump}) given an initial
condition $X_0 \in\mathbb{R}^d$, is well known under the above
assumptions, and for any $0 \leq T < \infty$, we have the
standard estimate
%
%
\begin{equation}\label{Xest}
\mathbf{E} \Bigl[{\sup_{0\leq t\leq T}}
|X_t|^2 \Bigr] < \infty.
\end{equation}

In what follows, we fix a finite time duration $[0,T]$. Let us
introduce some additional notation. We denote by:
\begin{itemize}
\item$\bolds{\mathcal{S}}^{\mathbf{2}}$\vspace*{0.5pt} the set of real-valued c\`
{a}dl\`{a}g adapted
processes $Y =
(Y_t)_{0\leq t\leq T}$ such that $ \|Y\|_{{\bolds{\mathcal{S}}^{\mathbf{2}}}}:=
(\mathbf{E} [ {\sup_{0\leq t\leq T}} |Y_t|^2 ]
)^{1/2} < \infty$.\vspace*{0.5pt}

\item$\mathbf{L}^{\mathbf{p}}(\mathbf{0},\mathbf{T})$, $p \geq1$, the
set of real-valued
processes $(\phi_t)_{0\leq t\leq T}$ such that\break $\mathbf{E} [\int_0^T
|\phi_t|^p \, dt ] < \infty$; and $\mathbf{L}^{\mathbf{p}}_{\mathbb
{F}}(\mathbf{0},\mathbf{T})$ is the
subset of $\mathbf{L}^{\mathbf{p}}(\mathbf{0},\mathbf{T})$ consisting
of adapted processes.

\item$\mathbf{L}^{\mathbf{p}}(\mathbf{W})$, $p \geq1$, the set of
$\mathbb
{R}^d$-valued $\mathcal{P}
$-measurable processes
$Z=(Z_t)_{0\leq t\leq T}$ such that $\|Z\|_{\mathbf{L}^{\mathbf
{p}}(\mathbf{W})} :=
(\mathbf{E} [ \int_0^T |Z_t|^p \,dt ] )^{1/p}
<\infty$.

\item$\mathbf{L}^{\mathbf{p}}(\tilde\mu)$, $p \geq1$, the set of
$\mathcal{P}\otimes\mathcal{E}$-measurable maps $U\dvtx\Omega\times
[0,T]\times
E\rightarrow\mathbb{R}$ such that $\| U\|_{{\mathbf{L}^{\mathbf
{p}}(\tilde\mu)}} :=
(\mathbf{E} [ \int_0^T\int_E |U_t(e)|^p \lambda(de) \, dt
] )^{1/p} < \infty$.

\item$\mathbf{A}^{\mathbf{2}}$ the closed subset of $\bolds{\mathcal
{S}}^{\mathbf{2}}$
consisting of
nondecreasing processes $K = (K_t)_{0\leq t\leq T}$ with $K_0
= 0$.
\end{itemize}

We are given four objects: (i) a terminal function, which is a
measurable function $g\dvtx\mathbb{R}^d \mapsto\mathbb{R}$ satisfying a
growth sublinear
condition
%
%
\begin{equation}\label{ggrowth} \sup_{x\in\mathbb{R}^d}
\frac{|g(x)|}{1+|x|} < \infty;
\end{equation}

(ii) a generator
function $f$, which is a measurable function
$f\dvtx\mathbb{R}^d\times\mathbb{R}\times\mathbb{R}^d\rightarrow
\mathbb{R}$ satisfying a growth
sublinear condition
%
%
\begin{equation}
\label{fgrowth}
\sup_{(x,y,z)\in\mathbb{R}^d\times\mathbb{R}\times\mathbb{R}^d}
\frac{|f(x,y,z)|}{1+|x|+|y|+|z|} < \infty
\end{equation}
and a\vspace*{1pt} uniform Lipschitz condition on $(y,z)$, that is, there exists a
constant $k_f$ such that for all $x \in\mathbb{R}^d$, $y,y' \in
\mathbb{R}$, $z,z'
\in\mathbb{R}^d$,
%
%
\begin{equation}
\label{flip}
|f(x,y,z)-f(x,y',z')| \leq k_f (|y-y'|+|z-z'|);
\end{equation}

(iii) a cost function, which is a measurable function $c\dvtx
\mathbb{R}^d\times\mathbb{R}\times\mathbb{R}^d\times E\rightarrow
\mathbb{R}$ satisfying a growth
sublinear condition
%
%
\begin{equation}
\label{cgrowth}
\sup_{(x,y,z,e)\in\mathbb{R}^d\times\mathbb{R}\times\mathbb
{R}^d\times E}
\frac{|c(x,y,z,e)|}{1+|x|+|y|+|z|} < \infty
\end{equation}
and a uniform Lipschitz condition on $(y,z)$, that is, there exists a
constant $k_c$ such that for all $x \in\mathbb{R}^d$, $y,y' \in
\mathbb{R}$, $z,z'
\in\mathbb{R}^d$, $e \in E$,
%
%
\begin{equation}
\label{clip}
|c(x,y,z,e)-c(x,y',z',e)| \leq k_c (|y-y'|+|z-z'|);
\end{equation}

(iv) a constraint function, which is a measurable map $h \dvtx
\mathbb{R}\times E \rightarrow\mathbb{R}$ s.t. for all $e \in
E$,
%
%
\begin{equation}
\label{hdec}
u \longmapsto h(u,e) \qquad\mbox{is nonincreasing},
\end{equation}
satisfying a Lipschitz condition on $u$, that is, there exists a
constant $k_h$ such that for all $u,u' \in\mathbb{R}$, $e \in
E$,
%
%
\begin{equation}
\label{hlip}
|h(u,e)-h(u',e)| \leq k_h |u-u'|
\end{equation}
and such that $\int_{E}|h(0,e)|\lambda(de) < +\infty$.

Let us now introduce our BSDE with constrained jumps: find a
quadruple $(Y,Z,U,K) \in\bolds{\mathcal{S}}^{\mathbf{2}}\times\mathbf
{L}^{\mathbf{2}}(\mathbf{W})\times
\mathbf{L}^{\mathbf{2}}(\tilde\mu)\times\mathbf{A}^{\mathbf{2}}$ satisfying
%
%
\begin{eqnarray}\label{BSDEgen}
\hspace*{31pt}
Y_t &=& g(X_T) + \int_t^T f(X_s,Y_s,Z_s) \,ds + K_T - K_t - \int
_t^T \langle Z_s, dW_s\rangle\nonumber\\[-8pt]\\[-8pt]
& & - \int_t^T\int_E \bigl( U_s(e)
- c(X_{s^-},Y_{s^-},Z_{s},e) \bigr) \mu(ds,de),\qquad 0 \leq t
\leq T, \mbox{ a.s.,}\nonumber
\end{eqnarray}
with
%
%
\begin{equation}\label{hcons}
h(U_t(e),e) \geq0,\qquad d\mathbf{P}\otimes dt\otimes\lambda(de),\mbox{
a.e.,}
\end{equation}
and such that for any other quadruple $(\tilde Y,\tilde
Z,\tilde U,\tilde K) \in\bolds{\mathcal{S}}^{\mathbf{2}}\times
\mathbf{L}^{\mathbf{2}}(\mathbf{W})\times
\mathbf{L}^{\mathbf{2}}(\tilde\mu)\times\mathbf{A}^{\mathbf{2}}$ satisfying
(\ref{BSDEgen}) and (\ref{hcons}), we have
\[
Y_t \leq\tilde Y_t,\qquad
0 \leq t \leq T, \mbox{ a.s.}
\]
We say that $Y$ is the minimal
solution to (\ref{BSDEgen}) and (\ref{hcons}). In the formulation of
Peng \cite{pen99}, one may sometimes say that $Y$ is the smallest supersolution
to (\ref{BSDEgen}) and (\ref{hcons}). We shall also say that $(Y,Z,U,K)$
is a minimal solution to (\ref{BSDEgen}) and~(\ref{hcons}), and we discuss
later the uniqueness of such quadruple.
\begin{Remark}
Since we are originally motivated by probabilistic
representation of QVIs, we put the BSDE with constrained jumps in
a Markovian framework. But all the results of Section \ref{sec3} about the
existence and approximation of a minimal solution hold true in a
general non-Markovian framework with the following standard
modifications: the terminal condition $g(X_T)$ is replaced by a
square integrable random variable $\xi\in
\mathbf{L}^{\mathbf{2}}(\bolds\Omega,\bolds{\mathcal{F}}_{\mathbf
{T}})$, the generator is a map $f$ from $\Omega
\times
[0,T] \times\mathbb{R}\times\mathbb{R}^d$ into $\mathbb{R}$,
satisfying a uniform Lipschitz
condition in $(y,z)$, and $f(\cdot,y,z) \in\mathbf{L}^{\mathbf
{2}}_\mathbb
{F}(\mathbf{0},\mathbf{T})$ for
all $(y,z) \in\mathbb{R}\times\mathbb{R}^d$, and the cost
coefficient is a map
$c$ from $\Omega\times[0,T] \times\mathbb{R}\times\mathbb
{R}^d\times E$ into $\mathbb{R}$,
satisfying a uniform Lipschitz condition in $(y,z)$, and
$c(\cdot,y,z,e) \in\mathbf{L}^{\mathbf{2}}_\mathbb{F}(\mathbf
{0},\mathbf{T})$ for all
$(y,z,e) \in\mathbb{R}\times\mathbb{R}^d\times E$.
\end{Remark}
\begin{Remark}
Without the $h$-constraint condition (\ref{hcons}) on jumps, we
have existence and uniqueness of a solution $(Y,Z,U,K)$ with $K = 0$ to
(\ref{BSDEgen}), from results on BSDE with jumps in
\cite{tanli94} and \cite{barbucpar97}. Here, under (\ref{hcons})
on jumps, it is not possible in general to have equality in
(\ref{BSDEgen}) with $K = 0$, and as usual in the BSDE literature
with constraint, we consider a nondecreasing process $K$ to have
more freedom. The problem is then to find a minimal solution to this
constrained BSDE, and the nondecreasing condition (\ref{hdec}) on $h$
is crucial for stating comparison principles needed in the
penalization approach. The primary example of constraint function
is $h(u,e) = -u$, that is, nonpositive jumps constraint, which is
actually equivalent to consider minimal solution to BSDE
(\ref{BSDE1}), as showed later.
\end{Remark}

\subsection{The case of nonpositive jump constraint} \label{paranonneg}

Let us recall the BSDE defined in the \hyperref[sec1]{Introduction}: find a triplet
$(Y,Z,K)\in\bolds{\mathcal{S}}^{\mathbf{2}}\times\mathbf{L}^{\mathbf
{2}}(\mathbf{W})\times\mathbf{A}^{\mathbf{2}}$
such that
%
%
\begin{eqnarray}\label{BSDEneg}
Y_t &=& g(X_T) + \int_t^T f(X_s,Y_s,Z_s)\, ds + K_T - K_t - \int_t^T
\langle Z_s, dW_s\rangle\nonumber\\[-8pt]\\[-8pt]
&&{} + \int_t^T\int_E
c(X_{s^-},Y_{s^-},Z_{s},e) \mu(ds,de),\qquad 0 \leq t \leq
T, \mbox{ a.s.,}\nonumber
\end{eqnarray}
such that for any other triplet $(\tilde
Y,\tilde Z,\tilde K) \in\bolds{\mathcal{S}}^{\mathbf{2}}\times\mathbf
{L}^{\mathbf{2}}(\mathbf{W})\times
\mathbf{A}^{\mathbf{2}}$ satisfying (\ref{BSDEneg}), it holds that
\[
Y_t \leq
\tilde Y_t,\qquad 0 \leq t \leq T, \mbox{ a.s.}
\]
We will call such
$Y$ [and, by a slight abuse of notation, $(Y,Z,K)$] the \textit{minimal
solution} to (\ref{BSDEneg}). We claim that this problem is actually
equivalent to problem (\ref{BSDEgen}) and (\ref{hcons}) in the case
$h(u,e) = -u$, corresponding to nonpositive jump constraint
condition
%
%
\begin{equation}\label{negjump}
U_t(e) \leq0,\qquad d\mathbf{P}\otimes dt\otimes\lambda(de),\mbox{ a.e.}
\end{equation}
Indeed,
let $(Y,Z,U,K)$ be any solution of
(\ref{BSDEgen}) and (\ref{negjump}). Define a process $\bar K$ by
$d\bar K_t = dK_t - \int_E U_s(e) \mu(dt,de)$, $0 \leq t\leq T$,
then $\bar K$ is nondecreasing, and the triplet $(Y,Z,\bar K)$
satisfies (\ref{BSDEneg}). It follows that
the minimal solution to (\ref{BSDEneg}) is smaller than the minimal
solution to (\ref{BSDEgen}) and (\ref{negjump}). We shall see in the
next section, by using comparison principles and penalization
approach, that equality holds, that is,
\[
\mbox{minimal solution $Y$
to } (\ref{BSDEneg}) = \mbox{minimal solution $Y$ to }
(\ref{BSDEgen}), (\ref{negjump}).
\]

We shall illustrate this result by considering a special case: when the
functions $f$ and $c$ do not depend on $y,z$ (i.e., the impulse
control case). In this case, one can obtain directly the existence
of a minimal solution to (\ref{BSDEgen})--(\ref{negjump}) and
(\ref{BSDEneg}) by duality methods involving the following set of
probability measures. Let $\mathcal{V}$ be the set of
$\mathcal{P}\otimes\mathcal{E}$-measurable essentially bounded
processes valued in
$(0,\infty)$, and given $\nu\in\mathcal{V}$, consider the
probability measure $\mathbf{P}^\nu$ equivalent to $\mathbf{P}$ on
$(\Omega,\mathcal{F}_T)$ with Radon--Nikodym density
%
%
\begin{equation}\label{defPnu}
\frac{d\mathbf{P}^\nu}{d\mathbf{P}} = \mathcal{E}_T \biggl( \int
_0^{\cdot} \int_E \bigl(\nu_t(e)-1\bigr)
\tilde\mu(dt,de) \biggr),
\end{equation}
where $\mathcal{E}_t(\cdot)$ is the Dol\'{e}ans--Dade
exponential. Notice that the Brownian motion $W$ remains a Brownian
motion under $\mathbf{P}^\nu$, which can then be interpreted as an
equivalent martingale measure for the ``asset'' price process $W$.
The effect of the probability measure $\mathbf{P}^\nu$, by Girsanov's
theorem, is to change the compensator $\lambda(de)\,dt$ of $\mu$ under
$\mathbf{P}$ to $\nu_t(e)\lambda(de)\,dt$ under $\mathbf{P}^\nu$.

In order to ensure that the problem is well defined, we need to assume:
\renewcommand{\theequation}{H1}
\begin{equation}\label{assuH1}\quad \ \ \
\mbox{There exists a triple $(\tilde Y,\tilde Z,\tilde K) \in\bolds{\mathcal
{S}} ^{\mathbf{2}}
\times\mathbf{L}^{\mathbf{2}}(\mathbf{W})\times\mathbf{A}^{\mathbf{2}}$
satisfying (\ref{BSDEneg}).}
\end{equation}
This assumption is standard and natural in the literature
on BSDE with constraints, and means equivalently here (when $f$
and $c$ do not depend on $y,z$) that one can find some constant
$\tilde y \in\mathbb{R}$, and $\tilde Z \in\mathbf{L}^{\mathbf
{2}}(\mathbf{W})$
such that
\[
\tilde y + \int_0^T \langle\tilde Z _s, dW_s\rangle\geq g(X_T)
+ \int_0^T f(X_s) \, ds + \int_0^T\int_E c(X_{s^-},e) \mu(ds,de),
\mbox{ a.s.}
\]
This equivalency can be proved by same arguments as
in \cite{cvikarson98}. Notice that assumption (\ref{assuH1}) may be not
satisfied as shown in Remark \ref{remexem}, in which case the
problem (\ref{BSDEneg}) is ill posed.
\begin{Theorem} \label{theonegjump}
Suppose that $f$ and $c$ do not depend on $y,z$ and (\ref{assuH1})
holds. Then, there exists a unique minimal solution $(Y,Z,K,U) \in
\bolds{\mathcal{S}}^{\mathbf{2}}\times\mathbf{L}^{\mathbf{2}}(\mathbf
{W})\times
\mathbf{L}^{\mathbf{2}}(\tilde\mu)\times\mathbf{A}^{\mathbf{2}}$, with
$K$ predictable, to
(\ref{BSDEgen})--(\ref{negjump}). Moreover, $(Y,Z,\bar K)$ is the
unique minimal solution to (\ref{BSDEneg}) with $\bar K_t = K_t - \int
_0^t\int_E U_s(e)\mu(ds,de)$, and $Y$ has the explicit
functional representation
\[
Y_t = \mathop{\operatorname{ess}\sup}_{\nu\in\mathcal{V}} \mathbf{E}^\nu\biggl[ g(X_T)
+ \int
_t^T f(X_s)\, ds +
\int_t^T \int_E c(X_{s^-},e) \mu(ds,de) \Big| \mathcal{F}_t
\biggr],
\]
for all $t\in[0,T]$.
\end{Theorem}
\begin{pf} First,\vspace*{1pt} observe that for any $(\tilde Y,\tilde Z,\tilde
U,\tilde K) \in\bolds{\mathcal{S}}^{\mathbf{2}}\times\mathbf
{L}^{\mathbf{2}}(\mathbf{W})\times\mathbf{L}^{\mathbf{2}}(\tilde\mu
)\times\mathbf{A}^{\mathbf{2}}$ [resp., $(\tilde Y,\tilde Z,\tilde
K) \in\bolds{\mathcal{S}}^{\mathbf{2}}\times\mathbf{L}^{\mathbf
{2}}(\mathbf{W})\times\mathbf{A}^{\mathbf{2}}$]
satisfying (\ref{BSDEgen})--(\ref{negjump}) [resp., (\ref{BSDEneg})], the
process
\[
\tilde Q_t := \tilde Y_t + \int_ 0^t f(X_s) \,ds +
\int_0^t\int_E c(X_{s^-},e) \mu(ds,de),\qquad 0 \leq t\leq T,
\]
is a $\mathbf{P}^\nu$-supermartingale, for all $\nu\in\mathcal{V}$,
where the
probability measure $\mathbf{P}^\nu$ was defined in (\ref{defPnu}). Indeed,
from (\ref{BSDEgen})--(\ref{negjump}) [resp., (\ref{BSDEneg})], we have
\begin{eqnarray*}
\tilde Q_t &=& \tilde Q_0 + \int_0^t
\langle
\tilde Z_s, dW_s\rangle-
\bar K_t\qquad \mbox{with } \bar K_t = \tilde K_t - \int_0^t U_s(e) \mu
(ds,de), \\
\biggl(\mbox{resp., } \tilde Q_t &=& \tilde
Q_0 + \int_0^t \langle
\tilde Z_s, dW_s\rangle- \tilde K_t\biggr),\qquad 0 \leq t\leq T.
\end{eqnarray*}
Now, by Girsanov's theorem, $W$ remains a Brownian motion under
$\mathbf{P}^\nu$, while from the boundedness of $\nu\in\mathcal{V}$, the
density $d\mathbf{P}^\nu/d\mathbf{P}$ lies in $L^2(\mathbf{P})$.
Hence, from
Cauchy--Schwarz inequality, the condition $\tilde Z \in\mathbf
{L}^{\mathbf{2}}(\mathbf{W})$, and Burkholder--Davis--Gundy inequality,
we get the
$\mathbf{P}^\nu$-martingale property of the stochastic integral $\int
\langle
\tilde Z, dW\rangle$, and so the $\mathbf{P}^\nu$-supermartingale
property of
$\tilde Q$ since $\bar K$ (resp., $\tilde K$) is nondecreasing. This
implies
\[
\tilde Y_t \geq\mathbf{E}^\nu\biggl[ \tilde Y_T +
\int_t^T f(X_s) \,ds + \int_t^T \int_E c(X_{s^-},e) \mu(ds,de) \Big|
\mathcal{F}_t \biggr]
\]
and thereby, from the arbitrariness of
$\mathbf{P}^\nu$, $\nu\in\mathcal{V}$, and since $\tilde
Y_T = g(X_T)$,
%
%
\setcounter{equation}{14}
\renewcommand{\theequation}{\arabic{section}.\arabic{equation}}
\begin{eqnarray}\label{YleqtildeY}\qquad
Y_t & := & \mathop{\operatorname{ess}\sup}_{\nu\in\mathcal{V}} \mathbf{E}^\nu
\biggl[ g(X_T) + \int_t^T f(X_s) \,ds
\nonumber\\[-8pt]\\[-8pt]
&&\hspace*{47.5pt}{}
+ \int_t^T \int_E c(X_{s^-},e)
\mu(ds,de) \Big| \mathcal{F}_t \biggr]
\leq\tilde Y_t.\nonumber
\end{eqnarray}

To show the converse, let us consider the process $Y$ defined in
(\ref{YleqtildeY}). By standard arguments as in \cite{elkque95}, the
process $Y$ can be considered in its c\`{a}dl\`{a}g modification, and
we also notice that $Y \in\bolds{\mathcal{S}}^{\mathbf{2}}$. Indeed,
by observing
that the choice of $\nu= 1$ corresponds to the probability
$\mathbf{P}^\nu= \mathbf{P}$, we have $\hat Y \leq Y \leq\tilde Y$,
where $(\tilde Y,\tilde Z,\tilde K) \in\bolds{\mathcal{S}}^{\mathbf
{2}}\times
\mathbf{L}^{\mathbf{2}}(\mathbf{W})\times\mathbf{A}^{\mathbf{2}}$ is a
solution to (\ref{BSDEneg}), and
\[
\hat Y_t = \mathbf{E} \biggl[ g(X_T) + \int_t^T f(X_s) \,ds + \int_t^T
\int_E c(X_{s^-},e) \mu(ds,de) \Big| \mathcal{F}_t \biggr].
\]
Thus, since $\hat Y$ lies in $\bolds{\mathcal{S}}^{\mathbf{2}}$ from
the linear growth
conditions on $g$, $f$ and $c$, and the estimate (\ref{Xest}), we
deduce that $Y \in\bolds{\mathcal{S}}^{\mathbf{2}}$. Now, by similar dynamic
programming arguments as in \cite{elkque95}, we see that the
process
%
%
\begin{equation}\label{defQ}\qquad
Q_t = Y_t + \int_ 0^t f(X_s) \,ds +
\int_0^t\int_E c(X_{s^-},e) \mu(ds,de),\qquad 0 \leq t\leq T,
\end{equation}
lies in $\bolds{\mathcal{S}}^{\mathbf{2}}$, and is a $\mathbf{P}^\nu
$-supermartingale, for all
$\nu\in\mathcal{V}$. Then, from the Doob--Meyer decomposition
of $Q$
under each $\mathbf{P}^\nu$, $\nu\in\mathcal{V}$, we obtain
%
%
\begin{equation}
\label{DoobQ}
Q_{t} = Y_{0}+M^{\nu}-K^{\nu},
\end{equation}
where $M^{\nu}$
is a $\mathbf{P}^\nu$-martingale, $M^\nu_{0} = 0$, and $K^{\nu
}$ is a
$\mathbf{P}^\nu$ nondecreasing predictable c\`{a}dl\`{a}g process with
$K^\nu_{0}=0$. Recalling that $W$ is a $\mathbf{P}^\nu$-Brownian motion,
and since $\tilde{\mu}^\nu(ds,de) := \mu(ds,de)-\nu_s(e)\lambda(de)\,ds$
is the compensated measure of
$\mu$ under~$\mathbf{P}^\nu$, the martingale representation theorem for
each $M^\nu$, $\nu\in\mathcal{V}$, gives the existence of predictable
processes $Z^{\nu}$ and $U^\nu$ such that
%
%
\begin{eqnarray}\label{decomQ}
Q_t &=& Y_0 + \int_0^t \langle Z_s^\nu,
dW_s\rangle\nonumber\\[-8pt]\\[-8pt]
& &{} +\int_0^t\int_E U_s^\nu(e)
\tilde\mu^\nu(ds,de) - K_t^\nu,\qquad 0 \leq t\leq T.\nonumber
\end{eqnarray}
By comparing the decomposition (\ref{decomQ}) under $\mathbf{P}^\nu$ and
$\mathbf{P}$ corresponding to $\nu= 1$, and identifying the
martingale parts and the predictable finite variation parts, we
obtain that $Z^\nu= Z^1 =: Z$, $U^\nu= U^1 =: U$
for all $\nu\in\mathcal{V}$, and
%
%
\begin{equation}\label{KnuK1}
K_t^\nu= K_t^1
- \int_0^t\int_E U_s(e)\bigl(\nu_s(e)-1\bigr)\lambda(de) \,ds,\qquad 0 \leq
t\leq T.
\end{equation}
Now, by writing the relation (\ref{decomQ}) with $\nu= \varepsilon>
0$, substituting the definition of $Q$ in
(\ref{defQ}), and since $Y_T = g(X_T)$, we obtain
%
%
\begin{eqnarray}\label{Yeps}
Y_t &=& g(X_T) + \int_t^T f(X_s) \,ds - \int_t^T \langle
Z_s,dW_s\rangle
\nonumber\\
& &{} - \int_t^T\int_E \bigl(U_s(e)-c(X_{s^-},e)\bigr)\mu(ds,de) \\
& &{} + \int_t^T\int_E U_s(e) \varepsilon\lambda(de)\,ds +
K_T^\varepsilon-K_t^\varepsilon,\qquad 0 \leq t\leq T. \nonumber
\end{eqnarray}
From (\ref{KnuK1}),
the process $K^\varepsilon$ has a limit as $\varepsilon$ goes to
zero, which is
equal to $K^0 = K^1+\int_0^{\cdot}\int_E U_s(e)\lambda(de)\,ds$, and
inherits from $K^\varepsilon$, the nondecreasing path and predictability
properties. Moreover, since $Q \in\bolds{\mathcal{S}}^{\mathbf{2}}$,
in the
decomposition (\ref{DoobQ}) of $Q$ under $\mathbf{P} = \mathbf
{P}^\nu$ for
$\nu= 1$, the process $M^1$ lies in $\bolds{\mathcal{S}}^{\mathbf{2}}$
and $K^1 \in\mathbf{A}^{\mathbf{2}}$. This implies that $Z \in\mathbf
{L}^{\mathbf{2}}(\mathbf{W})$, $U \in\mathbf{L}^{\mathbf{2}}(\tilde\mu
)$ and also that $K^0 \in\mathbf{A}^{\mathbf{2}}$.
By sending $\varepsilon$ to zero into~(\ref{Yeps}), we obtain that
$(Y,Z,U,K^0) \in\bolds{\mathcal{S}}^{\mathbf{2}}\times\mathbf
{L}^{\mathbf{2}}(\mathbf{W})\times
\mathbf{L}^{\mathbf{2}}(\tilde\mu)\times\mathbf{A}^{\mathbf{2}}$ is a
solution to (\ref{BSDEgen}). Let
us finally check that $U$ satisfies the constraint
%
%
\begin{equation}\label{unneg}
U_t(e) \leq0,\qquad d\mathbf
{P}\otimes
dt\otimes\lambda(de).
\end{equation}
We argue by contradiction by assuming
that the set $F=\{ (\omega,t,e)\in\Omega\times[0,T]\times E\dvtx
U_t(e) > 0\}$ has a strictly positive measure for $d\mathbf{P}\times
dt\times\lambda(de)$. For any $k > 0$, consider the process
$\nu_k = 1_{F^c} + (k+1)1_F$, which lies in $\mathcal{V}$. From
(\ref{KnuK1}), we have
\[
\mathbf{E}[K^{\nu_k}_{T}] = \mathbf{E}[K_{T}^1] -
k \mathbf{E} \biggl[ \int_0^T\int_E 1_F U_t(e) \lambda(de)\,dt \biggr]
< 0
\]
for $k$ large enough. This contradicts the fact that
$K^{\nu_k}_{T} \geq0$, and so (\ref{unneg}) is satisfied.
Therefore, $(Y,Z,U,K^0)$ is a solution to
(\ref{BSDEgen})--(\ref{negjump}), and it is a minimal solution from
(\ref{YleqtildeY}). $Y$ is unique by definition. The uniqueness of
$Z$ follows by identifying the Brownian parts and the finite
variation parts, and the uniqueness of $(U,K^0)$ is obtained by
identifying the predictable parts by recalling that the jumps of
$\mu$ are inaccessible. By denoting $\bar K^0 = K^0 -
\int_0^t\int_E U_s(e)\mu(ds,de)$, which lies in $\mathbf{A}^{\mathbf
{2}}$, we see
that $(Y,Z,\bar K^0)$ is a solution to (\ref{BSDEneg}), and it is
minimal by (\ref{YleqtildeY}). Uniqueness follows by identifying the
Brownian parts and the finite variation parts.
\end{pf}
\begin{Remark}
In Section \ref{secrelQVI}, we shall relate rigorously the
constrained BSDEs (\ref{BSDEgen}) and (\ref{hcons}) to QVIs. In
particular, the minimal solution $Y_t$ to
(\ref{BSDEgen})--(\ref{negjump}) or (\ref{BSDEneg}) is $Y_t=v(t,X_t)$
where $v$ is the value function of the impulse control problem
(\ref{impulse}). Together with the functional representation of $Y$
in Theorem \ref{theonegjump}, we then have the following relation at
time $t = 0$:
%
%
\setcounter{equation}{21}
\renewcommand{\theequation}{\arabic{section}.\arabic{equation}}
\begin{eqnarray}
\label{weakform}
v(0,X_0) &=& \sup_{\nu\in\mathcal{V}} \mathbf{E}^\nu\biggl[ g(X_T) +
\int_0^T f(X_s) \,ds \nonumber\\[-8pt]\\[-8pt]
&&\hspace*{32.76pt}{} +
\int_0^T \int_E c(X_{s^-},e) \mu(ds,de)\biggr].\nonumber
\end{eqnarray}
We then recover
a recent result obtained by Bouchard \cite{bou06}, who related
impulse controls to stochastic target problems in the case of a
finite set $E$. We may also interpret this result as follows.
Recall that the effect of the probability measure $\mathbf{P}^\nu$ is to
change the compensator $\lambda(de)\,dt$ of $\mu$ under $\mathbf{P}$ to
$\nu_t(e)\lambda(de)\,dt$ under~$\mathbf{P}^\nu$. Hence, by taking the
supremum over all $\mathbf{P}^\nu$, we formally expect to retrieve in
distribution law all the dynamics of the controlled process in
(\ref{controlX}) when varying
the impulse controls $\alpha$, which is confirmed by the equality
(\ref{weakform}).
\end{Remark}

Finally, we mention that the above duality and martingale methods
may be extended when the generator function $f$ is linear in $z$ by
using Girsanov's transformation.
Our main purpose is now to study the general case of
$h$-constraints on jumps, and nonlinear functions $f$ and $c$
depending on $y,z$.


\section{Existence and approximation by penalization}\label{sec3}

In this section, we prove the existence of a minimal solution to
(\ref{BSDEgen}) and (\ref{hcons}), based on approximation
via penalization. For each $n \in\mathbb{N}$, we introduce the penalized
BSDE with jumps
%
%
\setcounter{equation}{0}
\renewcommand{\theequation}{\arabic{section}.\arabic{equation}}
\begin{eqnarray}\label{BSDEpen}\qquad\quad
Y_t^n &=& g(X_T) + \int_t^T f(X_s,Y_s^n,Z_s^n) \,ds \nonumber\\
& &{} + n \int_t^T \int_E h^-(U^n_s(e),e)
\lambda(de)\,ds - \int_t^T \langle Z_s^n, dW_s\rangle\\
& &{} - \int_t^T\int_E \bigl(
U_s^n(e) - c(X_{s^-},Y_{s^-}^n,Z_{s}^n,e) \bigr) \mu(ds,de),\qquad
0 \leq t \leq T, \nonumber
\end{eqnarray}
where $h^-(u,e) = \max(-h(u,e),0)$ is the negative part of the function
$h$. Under
the Lipschitz and growth conditions on the coefficients $f$, $c$ and
$h$, we
know from the theory of BSDEs with jumps, see \cite{tanli94} and
\cite{barbucpar97}, that there exists a unique solution
$(Y^n,Z^n,U^n) \in\bolds{\mathcal{S}}^{\mathbf{2}}\times\mathbf
{L}^{\mathbf{2}}(\mathbf{W})\times
\mathbf{L}^{\mathbf{2}}(\tilde\mu)$ to (\ref{BSDEpen}). We define for
each $n \in\mathbb{N}$,
\[
K_t^n = n \int_0^t \int_E
h^-(U^n_s(e),e)\lambda(de)\,ds,\qquad 0 \leq t \leq T,
\]
which
is a nondecreasing process in $\mathbf{A}^{\mathbf{2}}$. The rest of
this section
is devoted to the convergence of the sequence $(Y^n,Z^n,U^n,K^n)_n$
to the minimal solution in which we are interested.

\subsection{Comparison results}\label{sec31}

We first state that the sequence $(Y^n)_n$ is nondecreasing. This
follows from a comparison theorem for BSDEs with jumps
whose generator is of the form $\tilde f(x,y,z,u) = f(x,y,z) + \int_E
\tilde h(u(e),e)\lambda(de)$ for some
nondecreasing function $\tilde h$, which covers our situation from the
nonincreasing condition on the constraint function $h$.
%
\begin{Lemma} \label{leminc}
The sequence $(Y^n)_n$ is nondecreasing, that is, for all $n \in\mathbb{N}
$, $Y^n_t \leq Y^{n+1}_t$, $0\leq t\leq T$, a.s.
\end{Lemma}
\begin{pf}
Define the sequence $(V^n)_{n}$ of $\mathcal{P}\otimes\mathcal{E}$-measurable
processes by
\[
V^n_{t}(e)  =  U^n_{t}(e)-c(X_{t^-},Y^n_{t^-},Z^n_{t},e),
\qquad
(t,e)\in(0,T]\times E,
\]
and
\[
V^n_{0}(e) =  U^n_{0}(e)-c(X_{0},Y^n_{0},Z^n_{0},e),\qquad e \in
E.
\]
From (\ref{BSDEpen}) and recalling that $X$ and $Y$ are
c\`{a}dl\`{a}g, we see that $(Y^{n},Z^n,V^n)$ is the unique solution
in $\bolds{\mathcal{S}}^{\mathbf{2}}\times\mathbf{L}^{\mathbf
{2}}(\mathbf{W})\times\mathbf{L}^{\mathbf{2}}(\tilde\mu)$ of the
BSDE with jumps:
\begin{eqnarray*}
Y_t^n &=& g(X_T) + \int_t^T
F_{n}(X_s,Y_s^n,Z_s^n,V_s^n) \,ds \\
& &{} - \int_t^T \langle Z_s^n, dW_s
\rangle- \int_t^T\int_E V^n_{s}(e) \tilde{\mu}(ds,de)
\end{eqnarray*}
with
$F_{n}(x,y,z,v) = f(x,y,z)+\int_{E} (nh^-(v(e)+c(x,y,z,e),e)-v(e)
)\lambda(de)$.
Since $h^-$ is nondecreasing, we have
%
\begin{eqnarray*}
&&F_{n}(t,x,y,z,v)-F_{n}(t,x,y,z,v') \\
&&\qquad =
\int_{E} \bigl\{\bigl(v'(e)-v(e)\bigr)+n\bigl[ h^-\bigl(v(e)+c(x,y,z,e),e\bigr)
\\
&&\qquad\quad\hspace*{99pt}{}- h^-\bigl(v'(e)+c(x,y,z,e),e\bigr)\bigr] \bigr\}\lambda(de) \\
&&\qquad \leq
\int_{E} \bigl\{ \bigl(-1+\mathbf{1}_{\{v(e)\geq v'(e)\}}
nk_h\bigr)\bigl(v(e)-v'(e)\bigr) \bigr\}\lambda(de).
\end{eqnarray*}

Moreover, since $F_{n+1}\geq F_{n}$, we can apply the comparison
Theorem 2.5 of \cite{roy06}, and obtain that
$Y^n_t \leq Y^{n+1}_t$, $0\leq t\leq T$, a.s.
\end{pf}

The next result shows that the sequence $(Y^n)_n$ is upper-bounded by
any solution to the constrained BSDE. Arguments in the proof involve
suitable change of probability measures $\mathbf{P}^\nu$, $\nu\in
\mathcal{V}
$, introduced in (\ref{defPnu}).
\begin{Lemma} \label{lemcompbor}
For any quadruple $(\tilde Y,\tilde Z,\tilde U,\tilde K) \in\bolds
{\mathcal{S}} ^{\mathbf{2}}\times\mathbf{L}^{\mathbf{2}}(\mathbf
{W})\times\mathbf{L}^{\mathbf{2}}(\tilde\mu)\times
\mathbf{A}^{\mathbf{2}}$ satisfying
(\ref{BSDEgen}) and (\ref{hcons}), and for all $n \in\mathbb{N}$,
we have
%
%
\begin{equation}\label{YnleqtildeY}
Y_t^n \leq\tilde Y_t,\qquad 0 \leq t\leq T, \mbox{ a.s.}
\end{equation}
Moreover, in the case: $h(u,e) = -u$, the inequality (\ref
{YnleqtildeY}) also holds for any triple
$(\tilde Y,\tilde Z,\tilde K) \in\bolds{\mathcal{S}}^{\mathbf{2}}\times
\mathbf{L}^{\mathbf{2}}(\mathbf{W})\times\mathbf{A}^{\mathbf{2}}$
satisfying (\ref{BSDEneg}).
\end{Lemma}
\begin{pf}
We state\vspace*{1pt} the proof for quadruple $(\tilde Y,\tilde Z,\tilde U,\tilde K)$
satisfying (\ref{BSDEgen}) and (\ref{hcons}). Same arguments are used in
the case: $h(u,e) = -u$ and $(\tilde Y,\tilde Z,\tilde K) \in\bolds
{\mathcal{S}}^{\mathbf{2}}\times\mathbf{L}^{\mathbf{2}}(\mathbf
{W})\times\mathbf{A}^{\mathbf{2}}$ satisfying
(\ref{BSDEneg}).

Denote $\bar{Y} = \tilde{Y}-Y^n$, $\bar{Z} = \tilde{Z}-Z^n$, $\bar{f} =
f(X,\tilde{Y},\tilde{Z})-f(X,Y^n,Z^n)$ and $\bar{c} = c(X_{.^-},\tilde
{Y}_{.^-},\tilde{Z},e)-c(X_{.^-},Y_{.^-}^n,Z^n,e)$.
Fix some $\nu\in\mathcal{V}$ (to be chosen later). We then have
%
\begin{eqnarray*}
\bar{Y}_t & = &
\int_t^T\bar{f}_{s}\,ds+\int_t^T \int_E\bar{c}_{s}\mu(ds,de)
-\int_t^T\langle\bar{Z}_{s},dW_{s}\rangle\\
& &{} -\int_t^T\int_E \{\tilde{U}_s(e)-U_s^n(e) \}\tilde{\mu
}^{\nu}(ds,de) 
\\
&&{} -\int_t^T\int_E \{\tilde{U}_s(e)-U_s^n(e) \}
\nu_{s}(e)\lambda(de)\,ds \\
& &{}
-n\int_{t}^T\int_{E}h^-(U^n_{s}(e),e)\lambda(de)\,ds
+\tilde{K}_T-\tilde{K}_t ,
\end{eqnarray*}
where $\tilde{\mu}^{\nu}(dt,de) = \mu(dt,de)-\nu_t(e)\lambda(de)\,dt$
denotes the compensated
measure of $\mu$ under $\mathbf{P}^\nu$. Let us then define the following
adapted processes:
\[
a_t = \frac{f(X_t,\tilde{Y}_t,\tilde
{Z}_t)-f(X_t,Y^n_t,\tilde{Z}_t)}{\bar{Y}_t}
\mathbf{1}_{\{\bar Y_t \neq0\}}
\]
and $b$ the $\mathbb{R}^d$-valued process defined by its $i$th
components, $i = 1,\ldots,d$:
\[
b_t^i =
\frac{f(X_t,Y^n_t,Z^{(i-1)}_t)-f(X_t,Y^n_t,Z^{(i)}_t)}{V_t^{i}}
\mathbf{1}_{\{V_t^{i}\neq0\}},
\]
where $Z_t^{(i)}$ is the $\mathbb{R}^d$-valued random vector whose
$i$ first components are those of $\tilde{Z}$ and whose $(d-i)$
lasts are those of $Z^n$, and $V_t^{i}$ is the $i$th component
of $Z^{(i-1)}_t-Z^{(i)}_t$. Let us also define the
$\mathcal{P}\otimes\mathcal{E}$-measurable processes $\delta$ in
$\mathbb{R}$ and $\ell$ in
$\mathbb{R}^d$ by
\[
\delta_t(e) =
\frac{c(X_{t^-},\tilde{Y}_{t^-},\tilde
{Z}_t,e)-c(X_{t^-},Y^n_{t^-},\tilde{Z}_t,e)}{\bar{Y}_r}
\mathbf{1}_{\{\bar{Y}_{t^-}\neq0\}}
\]
and
\[
\ell_r^i(e) =
\frac
{c(X_{t^-},Y^n_{t^-},Z^{(i-1)}_t,e)-c(X_{t^-},Y^n_{t^-},Z^{(i)}_t,e)}{V_t^{i}}
\mathbf{1}_{\{V_t^{i}\neq0\}}.
\]
Notice that the processes $a,b,\delta$ and $\ell$ are bounded by
the Lipschitz conditions on $f$ and $c$. Define also $\alpha_t^\nu=
a_t+\int_E \delta_t(e)\nu_t(e)\lambda(de)$, $\beta_t^\nu= b_t+\break\int
_E\ell_t(e)\nu_t(e)\lambda(de)$, which are bounded
processes since $a,b,\delta,\ell$ are bounded and $\lambda$ is a
finite measure on $E$, and denote $V_t^n(e) = \tilde{U}_t(e) -
U_t^n(e) - \delta_t(e)\bar{Y}_t - \ell_t(e)\cdot\bar{Z}_t$. With this
notation, and recalling that $h^-(\tilde U_s(e),e) = 0$
from the constraint condition (\ref{hcons}), we rewrite the BSDE for
$\bar Y$ as
%
\begin{eqnarray*}
\bar{Y}_t &=&\int_t^T
(\alpha_s^\nu\bar{Y}_s+\beta_{s^.}^\nu\bar{Z}_s )\,ds -
\int_t^T\langle\bar{Z}_s, dW_s\rangle
\\
&&{}
- \int_t^T\int_EV_s^n(e)\tilde{\mu}^{\nu}(ds,de)
+\tilde{K}_T-\tilde{K}_t \\
& &{}+
\int_t^T\int_E \{n[h^-(\tilde
{U}_s(e),e)-h^-(U_s^n(e),e)]\\
&&\hspace*{44.6pt}\hspace*{36.04pt}{} - \nu_s(e)[\tilde{U}_s(e)-U_s^n(e)]
\}\lambda(de)\,ds.
\end{eqnarray*}
Consider now the positive process $\Gamma^\nu$ solution to
the s.d.e.
%
%
\begin{equation}\label{dyngam}
d\Gamma_t^\nu= \Gamma_t^\nu(\alpha_t^\nu
dt+\langle\beta_t^\nu, dW_t\rangle),\qquad \Gamma_0^\nu
= 1,
\end{equation}
and notice that $\Gamma^\nu$ lies in $\mathcal{S}^2$ from the
boundeness condition on $\alpha^\nu$ and $\beta^\nu$. By It\^{o}'s
formula, we have
\begin{eqnarray*}
d\Gamma_t^\nu\bar{Y}_t& =&-\Gamma
_{t}^\nu\int_{E} \{n[h^-(\tilde
{U}_t(e),e)-h^-(U_t^n(e),e)]\\
&&\hspace*{73pt}{}-\nu_t(e)[\tilde{U}_t(e)-U_t^n(e)]
\}\lambda(de)\,ds \\
& &{}- \Gamma_t^\nu \,d\tilde{K}_t+\Gamma_t^\nu\langle
\bar
{Z}_t, dW_t\rangle
+\Gamma_{t}^\nu\bar{Y}_{t-}\langle\beta_{t}, dW_{t}\rangle
+\Gamma_t^\nu\int_EV_t^n(e)\tilde{\mu}^{\nu}(dt,de),
\end{eqnarray*}
which
shows that the process
\begin{eqnarray*}
&&\Gamma_t^{\nu} \bar{Y}_t + \int_0^t\Gamma_s^{\nu}
\int_{E}\{n[h^-(\tilde
{U}_s(e),e)-h^-(U_s^n(e),e)]\\
&&\hspace*{113.73pt}{} -\nu_s(e)[\tilde{U}_s(e)-
U_s^n(e)] \}\lambda(de)\,ds
\end{eqnarray*}
is a $\mathbf{P}^{\nu}$-supermartingale 
and so
\begin{eqnarray*}
\Gamma_t^{\nu} \bar{Y}_t &
\geq& \mathbf{E}^{\nu} \biggl[
\int_t^T\Gamma_s^{\nu}\int_{E} \{
n[h^-(\tilde{U}_s(e),e)-h^-(U_s^n(e),e)]\\
& &\hspace*{98.17pt}{} -\nu_s^\varepsilon(e)[\tilde{U}_s(e)-U_s^n(e)] \}
\lambda(de)\,ds \Big|\mathcal{F}_t \biggr].
\end{eqnarray*}
Now, from the Lipschitz condition on $h$, we see that the process $\nu
^\varepsilon$ defined by
\[
\nu^\varepsilon_t(e) =
\cases{\dfrac{n[h^-(\tilde{U}_t(e),e)-h^-(U_t^n(e),e)]}{\tilde {U}_t(e)-U_t^n(e)}, &\quad if $U_t^n(e)>\tilde{U}_{t}(e)$,\cr
&\quad and $h^-(U_t^n(e),e)>0$,\vspace*{2pt}\cr
\varepsilon, &\quad else,}
\]
is bounded and so lies in $\mathcal{V}$, and therefore by taking $\nu
= \nu^\varepsilon$, we obtain
%
%
\begin{eqnarray} \label{SMG}\hspace*{33pt}
\Gamma_t^{\nu^\varepsilon} \bar{Y}_t &\geq& -\varepsilon\mathbf
{E}^{\nu
^\varepsilon} \biggl[
\int_t^T\Gamma_s^{\nu^\varepsilon} \int_{E}[\tilde{U}_s(e)-U_s^n(e)]
\nonumber\\
& &\hspace*{80.1pt}{}\times \mathbf{1}_{\{\tilde{U}_s(e)\geq U_s^n(e)\}\cup\{
h^{-}(U_{s}^{n}(e),e)=0\}}\lambda(de)\,ds
\Big|\mathcal{F}_t \biggr] \\
& =: & -\varepsilon R_t^\varepsilon,\qquad 0 \leq t \leq T.\nonumber
\end{eqnarray}
From the conditional Cauchy--Schwarz inequality and Bayes formula we
have for all $t \in[0,T]$, $\varepsilon> 0$,
\begin{eqnarray*}
| R_t^\varepsilon| & \leq&\sqrt{
\mathbf{E} \biggl[\frac{Z_T^{\nu^\varepsilon}}{Z_t^{\nu
^\varepsilon
}} \int_t^T|\Gamma_s^{\nu^\varepsilon}|^2\,ds \Big|\mathcal
{F}_t \biggr]}
\\
& &\times{} \biggl( \mathbf{E} \biggl[\frac{Z_T^{\nu^\varepsilon
}}{Z_t^{\nu^\varepsilon}} \int_t^T \biggl(\int_{E}[\tilde
{U}_s(e)-U_s^n(e)]\\
&&\hspace*{81.5pt}{}\times \mathbf{1}_{\{\tilde{U}_s(e)\geq U_s^n(e)\}\cup\{
h^{-}(U_{s}^{n}(e),e)=0\}}\lambda(de) \biggr)^2\,ds
\Big|\mathcal{F}_t \biggr]\biggr)^{1/2} \\
&=: & R_t^{1,\varepsilon} R_t^{2,\varepsilon}.
\end{eqnarray*}
By definition of $\nu^\varepsilon$, we have for $\varepsilon\leq
n k_{h}$
\[
\frac{Z_T^{\nu^\varepsilon}}{Z_t^{\nu^\varepsilon}} \leq
\frac{Z_T^n}{Z_t^n}\exp\biggl(\int_t^T\int_{E}nk_{h}\lambda
(de)\,ds \biggr),
\]
where $Z^n$ is the solution to $dZ^n_t = Z^n_{t-}\int_E
(nk_{h}-1 )\tilde{\mu}(dt,de)$, $Z^n_0 = 1$. It follows
that for all
$t \in[0,T]$, $(R_t^{2,\varepsilon})_\varepsilon$ is
uniformly bounded for
$\varepsilon$ in a neighborhood of $0^+$. Similarly, using also the
boundedness of the coefficients $\alpha^{\nu^\varepsilon}$ and
$\beta^{\nu^\varepsilon}$ in the dynamics (\ref{dyngam}) of $\Gamma
^{\nu,\varepsilon}$, we deduce that
$(R_t^{1,\varepsilon})_\varepsilon$ and thus $(R_t^\varepsilon
)_\varepsilon$ is uniformly
bounded for $\varepsilon$ in a neighborhood of $0^+$.
Finally, since $\lim_{\varepsilon\rightarrow0} \Gamma_t^{\nu
^\varepsilon
} = \Gamma_t^{\nu^0} > 0$, by sending $\varepsilon$ to zero into
(\ref{SMG}), we conclude that $\bar{Y}_t \geq0$.
\end{pf}

\subsection{Convergence of the penalized BSDEs}\label{sec32}

We impose the following analogue of assumption (\ref{assuH1}):
\renewcommand{\theequation}{H2}
\begin{eqnarray}\label{assuH2}
\\[-19pt]
\begin{tabular}{p{327pt}}
There exists a quadruple $(\tilde
Y,\tilde Z,\tilde K,\tilde U) \in\bolds{\mathcal{S}}^{\mathbf{2}}\times
\mathbf{L}^{\mathbf{2}}(\mathbf{W})\times\mathbf{L}^{\mathbf{2}}(\tilde
\mu)\times
\mathbf{A}^{\mathbf{2}}$ satisfying (\ref{BSDEgen}) and (\ref{hcons}).
\end{tabular}\hspace*{-37pt}\nonumber
\end{eqnarray}

Assumption (\ref{assuH2}) ensures that the problem
(\ref{BSDEgen}) and (\ref{hcons}) is well posed. As indicated in Section
\ref
{paranonneg}, assumption (\ref{assuH2}) in the case $h(u,e) = -u$,
implies assumption (\ref{assuH1}). Since (\ref{assuH1}) is obviously stronger
than (\ref{assuH2}), these two assumptions are equivalent in the case
$h(u,e) = -u$.
We provide in Section \ref{secsuff} some discussion and sufficient
conditions under which (\ref{assuH2}) holds.
\begin{Remark} \label{remexem}
The following example shows that conditions (\ref{assuH1}) and
(\ref{assuH2}) may be not satisfied: consider the BSDEs
%
%
\setcounter{equation}{4}
\renewcommand{\theequation}{\arabic{section}.\arabic{equation}}
\begin{equation}
\label{BSDEexeneg}
Y_{t} = -\int_{t}^T\langle
Z_{s},dW_{s}\rangle+
\int_{t}^T\int_{E} c \mu(ds,de)+K_T-K_t
\end{equation}
and
%
%
\begin{equation}\label{BSDECE}\quad
\cases{Y_{t}=\displaystyle-\int_{t}^T\langle Z_{s},dW_{s}\rangle-\int_{t}^T\int
_{E}[U_{s}(e)-c]\mu(ds,de)+K_T-K_t,\vspace*{2pt}\cr
-U_{s}(e)\geq0,}
\end{equation}
where $c$ is a strictly positive constant, $c > 0$. Then, there
does not exist any solution to (\ref{BSDEexeneg}) or (\ref{BSDECE})
with component $Y \in\bolds{\mathcal{S}}^{\mathbf{2}}$. On the
contrary, we would
have
\[
Y_{0} \geq-\int_{0}^T\langle
Z_{s},dW_{s}\rangle+ c
\mu([0,T]\times E) \qquad\mbox{a.s.,}
\]
which implies
that for all $n \in\mathbb{N}^*$, $\nu\equiv n \in
\mathcal{V}$,
\[
Y_{0} \geq\mathbf{E}^\nu\biggl[-\int
_{0}^T\langle Z_{s},dW_{s}\rangle+c
\mu([0,T]\times E) \biggr] = c n\lambda(E)T.
\]
By sending
$n$ to infinity, we get the contradiction: $\|Y\|_{\mathcal{S}^2} =
\infty$.
\end{Remark}

We now establish a priori estimates, uniform on $n$, on the sequence
$(Y^n,Z^n,\break U^n$, $K^n)_n$.
\begin{Lemma} \label{lembor}
Under (\ref{assuH2}) [or (\ref{assuH1}) in the case: $h(u,e) = -u$], there
exists some constant $C$ such that
%
%
\setcounter{equation}{6}
\renewcommand{\theequation}{\arabic{section}.\arabic{equation}}
\begin{equation}\label{bounduni}
\|Y^n\|_{{\bolds{\mathcal{S}}^{\mathbf{2}}}} + \|Z^n\|_{{\mathbf
{L}^{\mathbf{2}}(\mathbf{W})}} + \|U^n\|
_{{\mathbf{L}^{\mathbf{2}}(\tilde\mu)}} +\|K^n\|_{{\bolds{\mathcal
{S}}^{\mathbf{2}}}} \leq C\qquad
\forall n \in\mathbb{N}.
\end{equation}
\end{Lemma}
\begin{pf}
In what follows, we shall denote $C>0$ to be a generic
constant depending only on $T$, the coefficients $f$, $c$, the
process $X$ and the bound for $\tilde Y$ in (\ref{assuH1}) or
(\ref{assuH2}), and which may vary from line to line.

Applying It\^{o}'s formula to $|Y_t^n|^2$, and observing that $K^n$ is
continuous and
$\Delta Y_t^n = \int_E \{U_s^n(e)-c(X_{s^-},Y^n_{s^-},Z_s^n,e)\}
\mu(\{t\},de)$, we have
\begin{eqnarray*}
\mathbf{E}|g(X_T)|^2 & = & \mathbf
{E}|Y_t^n|^2-2\mathbf{E}
\int_t^TY_s^nf(X_s,Y_s^n,Z_s^n)\,ds
\\ & &{} -2\mathbf{E}\int_t^TY_s^ndK^n_s+\mathbf{E}\int
_t^T|Z^n_s|^2\,ds \\
& &{} +
\mathbf{E}\int_t^T\int_E \{
|Y_{s^-}^n+U_s^n(e)\\
&&\hspace*{56.4pt}{}-c(X_{s^-},Y^n_{s^-},Z^n_{s},e)|^2-|Y_{s-}^n|^2
\}\lambda(de)\,ds.
\end{eqnarray*}
From the linear growth condition on $f$
and the inequality $Y_t^n \leq\tilde{Y}_t$ by Lemma
\ref{lemcompbor} under (\ref{assuH2}) [and also under (\ref{assuH1}) in the
case $h(u,e) = -u$], and using the inequality $2ab \leq\frac{1}{\alpha
}a^2 + \alpha b^2$ for any constant $\alpha>0$,
we have
\begin{eqnarray*}
& & \mathbf{E}|Y_t^n|^2+\mathbf{E}\int
_t^T|Z^n_s|^2\,ds+\mathbf{E}\int
_t^T\int_E|U_s^n(e)-
c(X_{s^-},Y_{s^-}^n,Z_{s}^n,e)|^2\lambda(de)\,ds \\
&&\qquad\leq
\mathbf{E}|g(X_T)|^2 + 2C\mathbf{E}\int_t^T|Y_s^n|
(1+|X_s|+|Y_s^n|+|Z_s^n| )\,ds \\
& &\qquad\quad{} - 2\mathbf{E}\int_{t}^T\int
_{E}Y_{s-}^{n}\bigl(U_s^n(e) -c(X_{s^-},Y_{s^-}^n,Z_{s}^n,e)\bigr)\lambda(de)\,ds
\\
& &\qquad\quad{} +\frac{1}{\alpha}\mathbf{E} \Bigl[{\sup_{t\in
[0,T]}}|\tilde{Y}_t|^2 \Bigr] + \alpha\mathbf{E}|K^n_T-K^n_t|^2.
\end{eqnarray*}
Using again the inequality $2ab \leq\frac{1}{\eta}a^2 + \eta b^2$ for
$\eta> 0$ yields
\begin{eqnarray*}
&& \mathbf{E}|Y_t^n|^2+\mathbf{E}\int_t^T|Z^n_s|^2\,ds
\\
&&\quad{}
+\frac{1-\eta}{2}\mathbf{E}\int_t^T\int_E|U_s^n(e)
-c(X_{s^-},Y_{s^-}^n,Z_{s}^n,e)|^2\lambda(de)\,ds \\
&&\qquad\leq \mathbf{E}|g(X_T)|^2 +2C\mathbf{E}\int_t^T|Y_s^n|
(1+|X_s|+|Y_s^n|+|Z_s^n| )\,ds \\
&&\qquad\quad{} + \frac{\lambda(E)}{\eta}\mathbf{E}\int_{t}^T|Y_{s}^{n}|^2\,ds
+\frac{1}{\alpha}\mathbf{E} \Bigl[{\sup_{t\in[0,T]}}|\tilde
{Y}_t|^2 \Bigr]
+\alpha\mathbf{E}|K^n_T-K^n_t|^2 \\
&&\qquad\leq C \biggl( 1 + \mathbf{E}\int_{t}^T |Y_{s}^{n}|^2 \,ds \biggr)
+ \frac{1}{2}\mathbf{E}\int_t^T |Z^n_s|^2 \,ds \\
&&\qquad\quad{} + \alpha\mathbf{E}|K^n_T-K^n_t|^2+\frac{\lambda(E)}{\eta
}\mathbf{E}\int
_{t}^T|Y_{s}^{n}|^2\,ds.
\end{eqnarray*}
Then, by using the inequality $(a-b)^2 \geq a^2/2 - b^2$, we get
%
%
\setcounter{equation}{7}
\renewcommand{\theequation}{\arabic{section}.\arabic{equation}}
\begin{eqnarray} \label{interYn}
&& \mathbf{E}|Y_t^n|^2+ \frac{1}{2} \mathbf{E}\int_t^T|Z^n_s|^2\,ds+
\frac{1-\eta
}{4} \mathbf{E}\int_t^T\int_E|U_s^n(e)|^2\lambda(de)\,ds \nonumber\\
&&\qquad\leq \frac{1-\eta}{2} \mathbf{E}\int_t^T\int
_E|c(X_{s^-},Y_{s^-}^n,Z_{s}^n,e)|^2\lambda(de)\,ds \nonumber\\
&&\qquad\quad{} + C \biggl( 1 + \mathbf{E}\int_{t}^T |Y_{s}^{n}|^2 \,ds \biggr) +
\alpha\mathbf{E}|K^n_T-K^n_t|^2 \nonumber\\[-8pt]\\[-8pt]
&&\qquad \leq C \biggl( 1 + \mathbf{E}\int_{t}^T |Y_{s}^{n}|^2 \,ds \biggr)\nonumber\\
&&\qquad\quad{} + C(1-\eta
)\mathbf{E}\int_t^T
|Z^n_s|^2\,ds\nonumber\\
&&\qquad\quad{} + \alpha\mathbf{E}|K^n_T-K^n_t|^2\nonumber
\end{eqnarray}
from the linear growth condition on $c$. Now, from the relation
\begin{eqnarray*}
K^n_T-K^n_t &=& Y_t^n-g(X_T)\\
&&{}-\int_t^T f(X_s,Y_s^n,Z_s^n)\,ds \\
&&{} + \int_t^T\int_E\bigl(U_s^n(e)-c(X_{s^-},Y_{s^-}^n,Z_{s}^n)\bigr)\mu
(ds,de)\\
&&{} + \int_t^T\langle Z^n_s, dW_s\rangle
\end{eqnarray*}
and the linear growth condition on $f$, $c$, there exists some positive
constant $C_1$ s.t.
%
%
\begin{eqnarray}\label{inegK}
&& \mathbf{E}|K^n_T-K^n_t|^2 \nonumber\\
&&\qquad \leq
C_1 \biggl(1+ \mathbf{E}|Y^n_t|^2+ \mathbf{E}\int_t^T (
|Y^n_s|^2+|Z_s^n|^2) \,ds\\
&&\hspace*{91.7pt}{} + \mathbf{E}\int_t^T \int_E
|U_s^n(e)|^2\lambda(de) \,ds \biggr). \nonumber
\end{eqnarray}
Hence, by choosing $\eta> 0$ s.t. $ (\frac{1}{2}-C(1-\eta
) )\wedge( \frac{1-\eta}{2} ) > 0$ and $\alpha> 0$ s.t.
$C_1\alpha< (\frac{1}{2}-C(1-\eta) )\wedge(
\frac{1-\eta}{2} )$, and plugging into (\ref{interYn}), we get
\begin{eqnarray*}
&&\mathbf{E}|Y_t^n|^2+\mathbf{E}\int_t^T|Z^n_s|^2\,ds+\mathbf{E}\int
_t^T\int
_E|U_s^n(e)|^2\lambda(de)\,ds\\
&&\qquad \leq C \biggl(1+\mathbf{E}\int_t^T|Y_s^n|^2
\,ds \biggr).
\end{eqnarray*}
By applying Gronwall's lemma to $t\mapsto\mathbf{E}|Y_t^n|^2$ and
(\ref{inegK}), we obtain
%
%
\begin{eqnarray}\label{ineg}
&&\sup_{0\leq t\leq T}\mathbf{E}|Y_t^n|^2+\mathbf{E}\int
_0^T|Z^n_s|^2\,ds \nonumber\\[-8pt]\\[-8pt]
&&\qquad{}+\mathbf{E}\int_0^T\int_E|U_s^n(e)|^2\lambda
(de)\,ds+\mathbf{E}
|K^n_T|^2 \leq C.\nonumber
\end{eqnarray}
Finally, by writing from (\ref{BSDEpen}) that
\begin{eqnarray*}
{\sup_{0\leq t\leq T}} |Y^n_{t}| & \leq& |g(X_{T})|+\int
_{0}^T|f(X_{s},Y_{s},Z_{s})|\,ds+K_{T}^n + \sup_{s\in[0,T]} \biggl|\int
_{0}^{T}\langle
Z_{s},dW_{s}\rangle\biggr|\\
& &{}
+\int_{0}^T\int_{E}|U_{s}^n(e)-c(X_{s^-},Y_{s^-},Z_{s},e)|\mu(ds,de),
\end{eqnarray*}
we obtain the required result from the Burkholder--Davis--Gundy
inequality, the linear growth condition on $f$, $c$ and
(\ref{ineg}).
\end{pf}
\begin{Remark}
A closer look at the proof leading to the estimate in
(\ref{bounduni}) shows that there exists a universal constant $C$,
depending only on $T$, and the linear growth condition constants of
$f$, $c$, such that for each $n \in\mathbb{N}$:
%
%
\setcounter{equation}{10}
\renewcommand{\theequation}{\arabic{section}.\arabic{equation}}
\begin{eqnarray}
\label{estiYngrowth}
\sup_{t\in[0,T]} \mathbf{E}[Y^n_t]^2 & \leq&
C \biggl( 1 + \mathbf{E}|g(X_T)|^2 \nonumber\\[-8pt]\\[-8pt]
&&\hspace*{13.6pt}{}+ \mathbf{E} \biggl[\int_0^T |X_t|^2 dt \biggr] +
\mathbf{E} \Bigl[{\sup_{t\in[0,T]}}|\tilde{Y}_t|^2 \Bigr] \biggr).
\nonumber
\end{eqnarray}
\end{Remark}
\begin{Lemma}\label{lemconvy}
Under (\ref{assuH2}) [or (\ref{assuH1}) in the case: $h(u,e) = -u$], the sequence of
processes $(Y_t^n)$ converges increasingly
to a process $(Y_t)$ with $Y \in\bolds{\mathcal{S}}^{\mathbf{2}}$. The
convergence also holds in $\mathbf{L}^{\mathbf{2}}_{\mathbb{F}}(\mathbf
{0},\mathbf{T})$ and for every
stopping time $\tau\in[0,T]$, the sequence of random variables
$(Y^n_\tau)$ converges to $Y_\tau$ in $\mathbf{L}^{\mathbf{2}}(\bolds
\Omega,\bolds{\mathcal{F}}_\tau)$, that is,
%
%
\setcounter{equation}{11}
\renewcommand{\theequation}{\arabic{section}.\arabic{equation}}
\begin{equation}\label{YnL2}
\lim_{n\rightarrow\infty} \mathbf{E} \biggl[\int_0^T |Y^n_t - Y_t|^2 \,dt
\biggr] = 0 \quad\mbox{and}\quad \lim_{n\rightarrow\infty} \mathbf
{E} [|Y^n_\tau- Y_\tau|^2 ] = 0.
\end{equation}
\end{Lemma}
\begin{pf}
From Lemmas \ref{leminc} and \ref{lemcompbor}, the
(nondecreasing) limit
%
%
\setcounter{equation}{12}
\renewcommand{\theequation}{\arabic{section}.\arabic{equation}}
\begin{equation}\label{limYn}
Y_t := \lim_{n\rightarrow\infty} Y_t^n,\qquad 0 \leq t \leq T,
\end{equation}
exists almost surely, and this defines an adapted process $Y$.
Moreover, by Lem\-ma~\ref{lembor} and convergence monotone theorem,
we have
\[
\mathbf{E} \Bigl[{\sup_{0\leq t\leq T}}|Y_t|^2 \Bigr]
<\infty.
\]
From the dominated convergence theorem, we also get
the convergence (\ref{YnL2}). It remains to check that the
process $Y$ has a c\`{a}dl\`{a}g modification. We first show that
$(Y^n)_n$ are quasi-martingales with uniformly bounded
conditional variations. That is, there exists a constant $C$ such
that, for any partition $\pi\dvtx0=t_0<t_1<\cdots< t_m=T$,
%
%
\begin{equation}
\label{quasimg}
\mathbf{E} \Biggl\{|Y^n_T| + \sum_{i=0}^{m-1}|\mathbf{E}\{
Y^n_{t_{i+1}}|\mathcal{F}_{t_i}\} -
Y^n_{t_i}| \Biggr\} \le C\qquad \forall\pi, \forall n.
\end{equation}
In fact, by (\ref{BSDEpen}), we have
%
\begin{eqnarray*}
&&\mathbf{E} \Biggl\{
\sum_{i=0}^{m-1} |\mathbf{E}\{Y^n_{t_{i+1}}|\mathcal{F}_{t_i}\} -
Y^n_{t_i}| \Biggr\}
\\
&&\qquad = \mathbf{E} \Biggl\{ \sum
_{i=0}^{m-1} \biggl|
\mathbf{E} \biggl[ \int_{t_i}^{t_{i+1}} f(X_s,Y_s^n,Z_s^n) \,ds \\
&&\qquad\quad\hspace*{45.3pt}{} + n \int_{t_i}^{t_{i+1}}\int_E
h^-(U^n_s(e),e)\lambda(de)\,ds \\
&&\qquad\quad{}\hspace*{45.3pt} -
\int_{t_i}^{t_{i+1}}\int_E \bigl( U_s^n(e) -
c(X_{s^-},Y_{s^-}^n,Z_{s}^n,e) \bigr)
\lambda(de)\,ds \Big| \mathcal{F}_{t_i} \biggr] \biggr| \Biggr\} \\
&&\qquad \le \mathbf{E} \biggl[ \int_0^T |f(X_s,Y_s^n,Z_s^n)|
\,ds \\
&&\qquad\quad\hspace*{10.5pt}{} + \int_0^T \int_E
|U_s^n(e) - c(X_{s^-},Y_{s^-}^n,Z_{s}^n,e)| \lambda(de)\,ds + K^n_T \biggr].
\end{eqnarray*}

Recall (\ref{ggrowth}), (\ref{fgrowth}) and (\ref{cgrowth}), we
have
%
\begin{eqnarray*}
&&\mathbf{E} \Biggl\{|Y^n_T| + \sum_{i=0}^{m-1}|\mathbf{E}\{
Y^n_{t_{i+1}}|\mathcal{F}_{t_i}\}
- Y^n_{t_i}| \Biggr\}
\\
&&\qquad\le C\mathbf{E} \biggl\{ 1+|X_T| + \int_0^T
[1+|X_s|+|Y_s^n|+|Z_s^n|] \,ds \\
&&\qquad\quad\hspace*{65.7pt}{} +\int_0^T \int_E |U_s^n(e)|\lambda
(de)\,ds + K^n_T \biggr\}.
\end{eqnarray*}
Applying (\ref{Xest}) and Lemma \ref{lembor}, we obtain
(\ref{quasimg}) immediately. Now by Meyer and Zheng \cite{mz} (or see
\cite{mz1}), there exists a subsequence $(Y^{n_k})_k$ and a
c\`{a}dl\`{a}g process $\tilde Y$ such that $(Y^{n_k})_k$ converges to
$\tilde Y$ in distribution. On the other hand, by (\ref{limYn}),
$(Y^{n_k})_k$ converges to $Y$, \textbf{P}-a.s. Then $Y$ and $\tilde Y$ have
the same distribution, and thus $Y$ is also c\`{a}dl\`{a}g.
\end{pf}

We now focus on the convergence of the diffusion and jump components
$(Z^n,U^n)$. In our context, we cannot prove the strong convergence
of $(Z^n,U^n)$ in $\mathbf{L}^{\mathbf{2}}(\mathbf{W})\times\mathbf
{L}^{\mathbf{2}}(\tilde\mu)$, and so
the strong convergence of $\int_0^t Z^n dW$ and $\int_0^t\int_E
U^n(s,e)\times\mu(ds,de)$ in $\mathbf{L}^{\mathbf{2}}(\bolds\Omega,\bolds
{\mathcal{F}}_t)$, see Remark
\ref{remstrong}. Instead, we follow and extend arguments of Peng
\cite{pen99}, and we shall prove that $(Z^n,U^n)$ converge in $\mathbf
{L}^{\mathbf{p}}(\mathbf{W})\times\mathbf{L}^{\mathbf{p}}(\tilde\mu)$,
for $1 \leq p < 2$.
First, we show the following weak convergence and decomposition
result.
\begin{Lemma}
Under (\ref{assuH2}) [or (\ref{assuH1}) in the case: $h(u,e) = -u$],
there exist $\phi\in\mathbf{L}^{\mathbf{2}}_{\mathbb{F}}(\mathbf
{0},\mathbf{T})$, $Z \in\mathbf{L}^{\mathbf{2}}(\mathbf{W})$, $V \in
\mathbf{L}^{\mathbf{2}}(\tilde\mu)$ and $K \in\mathbf{A}^{\mathbf{2}}$
predictable, such that the limit $Y$ in (\ref{limYn}) has the form
%
%
\setcounter{equation}{14}
\renewcommand{\theequation}{\arabic{section}.\arabic{equation}}
\begin{equation}\qquad
\label{decomY} Y_t = Y_0 - \int_0^t \phi_s \,ds - K_t +
\int_0^t\langle Z_s, dW_s\rangle+ \int_0^t \int_E V_s(e) \mu(ds,de)
\end{equation}
for all $t\in[0,T]$. Moreover, in the above decomposition of $Y$, the
components $Z$ and
$V$ are unique, and are, respectively, the weak limits of $(Z^n)$ in
$\mathbf{L}^{\mathbf{2}}(\tilde\mu)$ and of $(V^n)$ in $\mathbf
{L}^{\mathbf{2}}(\tilde\mu)$
where $V_t^n(e) = U_t^n(e) - c(X_{t^-},Y_{t^-}^n,Z_t^n,e)$,
$\phi$ is the weak limit in $\mathbf{L}^{\mathbf{2}}_\mathbb{F}(\mathbf
{0},\mathbf{T})$ of a
subsequence of
$(f^n) := (f(X,Y^n,Z^n))$, and $K$ is the weak limit in $\mathbf
{L}^{\mathbf{2}}_\mathbb{F}(\mathbf{0},\mathbf{T})$ of a subsequence of $(K^n)$.
\end{Lemma}
\begin{pf}
By Lemma \ref{lembor}, and the linear growth
conditions on $f$, $c$ together with (\ref{Xest}), the sequences
$(f^n)$, $(Z^n)$, $(V^n)$ are weakly compact, respectively, in $\mathbf
{L}_{\mathbb{F}}^{\mathbf{2}}(\mathbf{0},\mathbf{T})$, $\mathbf
{L}^{\mathbf{2}}(\mathbf{W})$ and $\mathbf{L}^{\mathbf{2}}(\tilde\mu)$.
Then, up
to a subsequence, $(f^n)$, $(Z^n)$, $(V^n)$ converge weakly to
$\phi$, $Z$ and $V$. By It\^{o} representation of martingales, we
then get the following weak convergence in $\mathbf{L}^{\mathbf
{2}}(\bolds\Omega,\bolds{\mathcal{F}}_\tau)$ for each stopping time
$\tau\leq T$:
\begin{eqnarray*}
\int_0^\tau f_s^n \,ds &\rightharpoonup& \int_0^\tau\phi_s \,ds,\qquad
 \int_0^\tau
\langle Z_s^n, dW_s\rangle\rightharpoonup\int_0^\tau
\langle Z_s,
dW_s\rangle, \\
\int_0^\tau\int_E V_s^n(e) \mu(ds,de) &\rightharpoonup&
\int_0^\tau\int_E V_s(e) \mu(ds,de).
\end{eqnarray*}
Since we have from (\ref{BSDEpen})
%
%
\setcounter{equation}{15}
\renewcommand{\theequation}{\arabic{section}.\arabic{equation}}
\begin{eqnarray}\label{dynforYn}
K_\tau^n &=& - Y_\tau^n +
Y_0^n - \int_0^\tau f_s^n \,ds \nonumber\\[-8pt]\\[-8pt]\
& &{} + \int_0^\tau\langle Z_s^n, dW_s\rangle+
\int_0^\tau\int_E V_s^n(e) \mu(ds,de),\nonumber
\end{eqnarray}
we also have the weak convergence in $\mathbf{L}^{\mathbf{2}}(\bolds
\Omega,\bolds{\mathcal{F}}_\tau)$:
%
%
\begin{eqnarray}
\label{defK} K_\tau^n \rightharpoonup K_\tau& := & -
Y_\tau+ Y_0 - \int_0^\tau\phi_s \,ds \nonumber\\[-8pt]\\[-8pt]
& &{} + \int_0^\tau\langle Z_s,
dW_s\rangle+ \int_0^\tau\int_E V_s(e) \mu(ds,de).\nonumber
\end{eqnarray}
The process $K$
inherits from $K^n$ the nondecreasing path property, is square
integrable, c\`{a}dl\`{a}g and adapted from (\ref{defK}), and so lies in
$\mathbf{A}^{\mathbf{2}}$. Moreover, by dominated convergence theorem,
we see that
$K^n$ converges weakly to $K$ in $\mathbf{L}^{\mathbf{2}}(\mathbf
{0},\mathbf{T})$. Since $K^n$ is
continuous, and so predictable, we deduce that $K$ is also
predictable, and we obtain the decomposition (\ref{decomY}) for $Y$.
The uniqueness of $Z$ follows by identifying the Brownian parts and
finite variation parts, and the uniqueness of $V$ is then obtained
by identifying the predictable parts and by recalling that the jumps
of $\mu$ are inaccessible. We conclude that $(Z,V)$ is uniquely
determined in (\ref{decomY}), and thus the whole sequence $(Z^n,V^n)$
converges weakly to $(Z,V)$ in $\mathbf{L}^{\mathbf{2}}(\mathbf
{W})\times\mathbf{L}^{\mathbf{2}}(\tilde\mu)$.
\end{pf}

The sequence $(U^n)$ is bounded in $\mathbf{L}^{\mathbf{2}}(\tilde\mu
)$, and so,
up to a subsequence, converges weakly to some $U \in\mathbf{L}^{\mathbf
{2}}(\tilde\mu)$. The next step is to show that the whole sequence
$(U^n)$ converges to $U$ and to identify in the decomposition
(\ref{decomY}) $\phi_t$ with $f(X_t,Y_t,Z_t)$, and $V_t(e)$ with
$U_t(e)-c(X_{t^-},Y_{t^-},Z_t,e)$. Since $f$ and $c$ are nonlinear,
we need a result of strong convergence for $(Z^n)$ and $(U^n)$ to
enable us to pass the limit in $f(X_t,Y_t^n,Z_t^n)$ as well as in
$U_t^n(e)-c(X_{t^-},Y_{t^-}^n,Z_t^n,e)$, and to eventually prove
the convergence of the penalized BSDEs to the minimal solution of
our jump-constrained BSDE. We shall borrow a useful technique of
Peng \cite{pen99} to carry out this task.
%
\begin{Theorem} \label{thmmain1}
Under (\ref{assuH2}), there exists a unique minimal solution
$(Y,Z,\break U,K) \in\bolds
{\mathcal{S}}^{\mathbf{2}}\times\mathbf{L}^{\mathbf{2}}(\mathbf
{W})\times\mathbf{L}^{\mathbf{2}}(\tilde
\mu)\times\mathbf{A}^{\mathbf{2}}$ with $K$ predictable, to (\ref
{BSDEgen}) and (\ref
{hcons}). $Y$ is the increasing limit of $(Y^n)$ in (\ref{limYn}) and
also in $\mathbf{L}_{\mathbb{F}}^{\mathbf{2}}(\mathbf{0},\mathbf{T})$,
$K$ is the weak limit of $(K^n)$ in $\mathbf{L}^{\mathbf{2}}_{\mathbb
{F}}(\mathbf{0},\mathbf{T})$, and
for any
$p \in[1,2)$,
\[
\|Z^n-Z\|_{{\mathbf{L}^{\mathbf{p}}(\mathbf{W})}} + \|U^n-U\|_{{\mathbf
{L}^{\mathbf{p}}(\tilde\mu)}}
\longrightarrow0
\]
as $n$ goes to infinity. Moreover, in the case: $h(u,e) = -u$,
$(Y,Z,\bar K)$ is the unique minimal solution to (\ref{BSDEneg}) with
$\bar K_t = K_t-\int_0^t\int_E U_s(e)\mu(ds,de)$, and this holds
true under (\ref{assuH1}). Consequently, the minimal solution $Y$ to
(\ref{BSDEneg}) and to (\ref{BSDEgen})--(\ref{negjump}) are the same.
\end{Theorem}
\begin{pf}
We apply It\^{o}'s formula to $|Y^n_{t}-Y_{t}|^2$ on a
subinterval $(\sigma,\tau]$, with $0 \leq\sigma< \tau\leq T$, two
stopping times. Recall the decomposition (\ref{decomY}), (\ref
{dynforYn}) of $Y$, $Y^n$, and observe that $K^n$ is continuous, and
$\Delta(Y_t^n-Y_t) = \Delta K_t + \int_E (V_t^n(e)-V_t(e))\mu
(\{t\},de)$. We then have
\begin{eqnarray*}
&&\mathbf{E}|Y^n_{\tau}-Y_{\tau}|^2 \\
&&\qquad = \mathbf{E}|Y^n_{\sigma
}-Y_{\sigma}|^2+
\mathbf{E}\int_{\sigma}^\tau|Z^n_{s}-Z_{s}|^2\,ds
+2\mathbf{E}\int_{\sigma}^\tau[Y^n_{s}-Y_{s}][\phi_{s}-f_s^n]\,ds\\
&&\qquad\quad{}
- 2\mathbf{E}\int_{\sigma}^\tau[Y^n_{s}-Y_{s}]\,dK^n_{s}
+2\mathbf{E}\int_{(\sigma,\tau]}[Y^n_{s^-}-Y_{s^-}]\,dK_{s}+\mathbf
{E}\sum_{t\in
(\sigma,\tau]}|\Delta K_{t}|^2 \\
&&\qquad\quad{}
+ \mathbf{E}\int_{(\sigma,\tau]}\int
_{E}[|Y^n_{s^-}-Y_{s^-}+V^n_{s}(e)-V_{s}(e)|^2-|Y^n_{s^-}-Y_{s^-}|^2]\mu
(ds,de) \\
&&\qquad = \mathbf{E}|Y^n_{\sigma}-Y_{\sigma}|^2+ \mathbf{E}\int
_{\sigma}^\tau
|Z^n_{s}-Z_{s}|^2\,ds+
2\mathbf{E}\int_{\sigma}^\tau[Y^n_{s}-Y_{s}][\phi_{s}-f_s^n]\,ds\\
&&\qquad\quad{} - 2\mathbf{E}\int_{\sigma}^\tau[Y^n_{s}-Y_{s}]\,dK^n_{s}
+ 2\mathbf{E}\int_{(\sigma,\tau]}[Y^n_{s^-}-Y_{s^-}+\Delta
K_{s}]\,dK_{s} \\
&&\qquad\quad{} - \mathbf{E}\sum_{t\in(\sigma,\tau]}|\Delta K_{t}|^2 +
\mathbf{E}\int
_{\sigma}^\tau\int_{E}|V^n_{s}(e)-V_{s}(e)|^2 \lambda(de)\,ds \\
&&\qquad\quad{} + 2 \mathbf{E}\int_{\sigma}^\tau\int_{E} (Y^n_{s}-Y_{s})
\bigl(V_s^n(e)-V_s(e)\bigr)\lambda(de)\,ds.
\end{eqnarray*}
Since $(Y^n_{s}-Y_{s})\,dK^n_{s} \leq0$, and by using the
inequality $2ab \geq-\frac{a^2}{2}-2b^2$ with $a =
V^n_{s}(e)-V_{s}(e)$ and $b = Y^n_{s}-Y_{s}$, we obtain
%
%
\setcounter{equation}{17}
\renewcommand{\theequation}{\arabic{section}.\arabic{equation}}
\begin{eqnarray}\label{estipeng}
&& \mathbf{E}\int_{\sigma}^\tau|Z^n_{s}-Z_{s}|^2\,ds+\frac
{1}{2}\mathbf{E}\int
_{\sigma}^\tau\int_{E} |V^n_{s}(e)-V_{s}(e)|^2\lambda(de)\,ds
\nonumber\\
&&\qquad \leq \mathbf{E}|Y^n_{\tau}-Y_{\tau}|^2 + 2 \mathbf{E}\int
_{\sigma}^\tau
|Y^n_{s}-Y_{s}|^2\,ds\nonumber\\[-8pt]\\[-8pt]
&&\qquad\quad{} + 2\mathbf{E}\int_{\sigma}^\tau
|Y^n_{s}-Y_{s}||\phi
_{s}-f_s^n|\,ds \nonumber\\
&&\qquad\quad{} + 2\mathbf{E}\int_{(\sigma,\tau]}
|Y^n_{s^-}-Y_{s^-}+\Delta
K_{s}|\,dK_{s} + \mathbf{E}\sum_{t\in(\sigma,\tau]}|\Delta
K_{t}|^2.\nonumber
\end{eqnarray}
The first two terms of the right-hand side of (\ref{estipeng}) converge
to zero by (\ref{YnL2}) in Lemma \ref{lemconvy}. The third term
also tends to zero since $(\phi-f^n)_n$ is bounded in $\mathbf
{L}^{\mathbf{2}}(\mathbf{0},\mathbf{T})$, and so by Cauchy--Schwarz inequality
%
%
\begin{equation}\quad
\label{Yncauchy}
\mathbf{E}\int_{0}^T |Y^n_{s}-Y_{s}||\phi_{s}-f_s^n|\,ds
\leq C \biggl(\mathbf{E}
\int_{0}^T
|Y^n_{s}-Y_{s}|^2\,ds \biggr)^{1/2} \rightarrow0.
\end{equation}
For the fourth term, we notice that the jumps of $Y^n$ are
inaccessible since they are determined by the Poisson random measure
$\mu$. Thus, the predictable projection of $Y^n$ is
$^pY^n_{t}=Y^n_{t^-}$. Similarly, from (\ref{decomY}), and since $K$
is predictable, we see that $^pY_{t}=Y_{t^-}-\Delta K_{t}$. Since
$Y^n$ increasingly converges to $Y$, then $^pY^n$ also
increasingly converges to $^pY$, and by the dominated convergence
theorem, we obtain
%
%
\begin{equation}\label{Ynsaut}
\lim_{n\rightarrow\infty}
\mathbf{E}\int_{(0,T]}|Y^n_{s^-}-Y_{s^-}+\Delta K_{s}|\,dK_{s} = 0.
\end{equation}
For the last term in (\ref{estipeng}), we apply Lemma 2.3 in
\cite{pen99} to the predictable nondecreasing process $K$: for any
$\delta,\varepsilon> 0$, there exists a finite number of pairs of
stopping times $(\sigma_k,\tau_k)$, $k = 0,\ldots,N$, with $0 < \sigma
_k \leq\tau_k \leq T$, such that all the
intervals $(\sigma_k,\tau_k]$ are disjoint and
%
%
\begin{equation}
\label{Keps}
\mathbf{E}\sum_{k=0}^N(\tau_k-\sigma_k) \geq T-\frac{\varepsilon
}{2},\qquad
\mathbf{E}\sum_{k=0}^N \sum_{\sigma_k<t\leq\tau_k} (\Delta K_t)^2
\leq
\frac{\varepsilon\delta}{3}.
\end{equation}
We should note that in \cite{pen99} the filtration is Brownian,
therefore it is continuous, and hence each
stopping time $\sigma_{k}$ can be approximated by a sequence of
announceable stopping times. In our case the stopping times
$\sigma_k$'s are constructed as the successive times of jumps of the
predictable process $K$ with size bigger than some given positive
level, the approximation of $\sigma_k$ by announceable stopping
times is again possible. We can thus argue exactly the same way as
in Lemma 2.3 in \cite{pen99} to derive both estimates in
(\ref{Keps}).

We now apply estimate (\ref{estipeng}) for each $\sigma= \sigma
_k$ and $\tau= \tau_k$,
and then take the sum over $k = 0,\ldots,N$. It follows that
\begin{eqnarray*}
\hspace*{-2pt}&& \sum_{k=0}^N \mathbf{E}\int_{\sigma_k}^{\tau_k} |Z^n_{s}-Z_{s}|^2\,ds
+\frac{1}{2} \sum_{k=0}^N \mathbf{E}\int_{\sigma_k}^{\tau_k} \int_{E}
|V^n_{s}(e)-V_{s}(e)|^2\lambda(de)\,ds \nonumber\\
\hspace*{-2pt}&&\qquad \leq \sum_{k=0}^N \mathbf{E}|Y^n_{\tau_k}-Y_{\tau_k}|^2 + 2
\mathbf{E}\int
_{0}^T |Y^n_{s}-Y_{s}|^2\,ds + 2\mathbf{E}\int_{0}^T
|Y^n_{s}-Y_{s}||\phi_{s}-f_s^n|\,ds
\nonumber\\
\hspace*{-2pt}&&\qquad\quad{} + 2\mathbf{E}\int_{(0,T]} |Y^n_{s^-}-Y_{s^-}+\Delta
K_{s}|\,dK_{s} +
\sum_{k=0}^N \mathbf{E}\sum_{t\in(\sigma_k,\tau_k]}|\Delta K_{t}|^2.
\end{eqnarray*}
From the convergence results in Lemma \ref{lemconvy}, (\ref{Yncauchy})
and (\ref{Ynsaut}),
we deduce that
\begin{eqnarray*}
&& \limsup_{n\rightarrow\infty}\sum_{k=0}^N \mathbf{E}\int
_{\sigma
_{k}}^{\tau_{k}} |Z^n_{s}-Z_{s}|^2\,ds
+\frac{1}{2}\sum_{k=0}^N\mathbf{E}\int_{\sigma_{k}}^{\tau
_{k}}\int_{E}
|V^n_{s}(e)-V_{s}(e)|^2\lambda(de)\,ds \\
&&\qquad \leq \sum_{k=0}^N\mathbf{E}\sum_{t\in(\sigma_{k},\tau
_{k}]}|\Delta
K_{t}|^2 \leq\frac{\varepsilon\delta}{3}.
\end{eqnarray*}
Thus, there exists an integer $\ell_{\varepsilon\delta}>0$ such
that for all $n\geq\ell_{\varepsilon\delta}$, we have
\[
\sum_{k=0}^N \mathbf{E}\int_{\sigma_{k}}^{\tau_{k}}
|Z^n_{s}-Z_{s}|^2\,ds+\frac{1}{2}\sum_{k=0}^N\mathbf{E}\int_{\sigma
_{k}}^{\tau_{k}}\int_{E}
|V^n_{s}(e)-V_{s}(e)|^2\lambda(de) \,ds \leq
\frac{\varepsilon\delta}{2}.
\]
This implies
\[
dt\otimes\mathbf{P} \Biggl[(s,\omega)\in
\bigcup_{k=0}^N(\sigma_{k}(\omega),\tau_{k}(\omega)] \times\Omega\dvtx
|Z^n_{s}(\omega)-Z_{s}(\omega)|^2 \geq\delta\Biggr] \leq\frac
{\varepsilon}{2}
\]
and
\begin{eqnarray*}
&&dt\otimes\lambda\otimes\mathbf{P} \Biggl[(s,e,\omega)\in\bigcup
_{k=0}^N(\sigma_{k}(\omega),\tau_{k}(\omega)] \\
&&\hspace*{53.5pt}\hspace*{62.3pt}{}\times\Omega\times E\dvtx
|V^n_{s}(e,\omega)-V_{s}(e,\omega)|^2\geq\delta\Biggr] \leq
\varepsilon.
\end{eqnarray*}
Together with (\ref{Keps}), it follows that
\[
dt\otimes\mathbf{P} \bigl[(s,\omega)\in[0,T]\times\Omega\dvtx
|Z^n_{s}(\omega)-Z_{s}(\omega)|^2\geq\delta\bigr] \leq
\varepsilon
\]
and
\begin{eqnarray*}
&&dt\otimes\lambda\times\mathbf{P} \bigl[(s,e,\omega)\in[0,T]\times E
\times\Omega\dvtx
|V^n_{s}(e,\omega)-V_{s}(e,\omega)|^2\geq\delta\bigr]\\
&&\qquad \leq
\varepsilon\bigl(1+\lambda(E)\bigr).
\end{eqnarray*}
We deduce that for all $\delta>0$
\[
\lim_{n\rightarrow\infty}
dt\otimes\mathbf{P} \bigl[(s,\omega)\in[0,T]\times\Omega\dvtx
|Z^n_{s}(\omega)-Z_{s}(\omega)|^2\geq\delta\bigr] = 0
\]
and
\[
\lim_{n\rightarrow\infty}dt\otimes\lambda\otimes\mathbf{P}
\bigl[(s,e,\omega)\in[0,T]\times
E \times\Omega\dvtx
|V^n_{s}(e,\omega)-V_{s}(e,\omega)|^2\geq\delta\bigr] = 0.
\]
This means that the sequences $(Z^n)_{n}$ and $(V^n)_{n}$
converge in measure, respectively, to $Z$ and $V$. Since they are
bounded, respectively, in $\mathbf{L}^{\mathbf{2}}(\mathbf{W})$ and
$\mathbf{L}^{\mathbf{2}}(\tilde\mu)$, they are uniformly integrable in
$\mathbf{L}^{\mathbf{p}}(\mathbf{W})$
and $\mathbf{L}^{\mathbf{p}}(\tilde\mu)$ for any $p \in[1,2)$,
respectively. Thus, $(Z^n)$ and $(V^n)$ converge strongly to $Z$
and $V$ in $\mathbf{L}^{\mathbf{p}}(\mathbf{W})$ and $\mathbf
{L}^{\mathbf{p}}(\tilde\mu)$,
respectively. Recalling that $U_t^n(e) = V_t^n(e) +
c(X_{t^-},Y_{t^-}^n,Z_t^n,e)$, and by the Lipschitz condition on
$c$, we deduce that the sequence $(U^n)$ converges strongly in
$\mathbf{L}^{\mathbf{p}}(\tilde\mu)$, for $p \in[1,2)$, to $U$ defined by
\[
U_{t}(e) = V_{t}(e) + c(X_{t^-},Y_{t^-},Z_{t},e),\qquad 0 \leq
t\leq T, e \in E.
\]
By the Lipschitz condition on $f$, we also have the strong convergence
in $\mathbf{L}^{\mathbf{p}}_\mathbb{F}(\mathbf{0},\mathbf{T})$ of
$(f^n) = (f(X,Y^n,Z^n))$ to $f(X,Y,Z)$. Since $\phi$ is the weak
limit of $(f^n)$ in $\mathbf{L}^{\mathbf{2}}_\mathbb{F}(\mathbf
{0},\mathbf{T})$, we deduce that
$\phi= f(X,Y,Z)$. Therefore, with the decomposition (\ref
{decomY}) and since $Y_T = \lim_n Y_T^n = g(X_T)$, we obtain
immediately
that $(Y,Z,U,K)$ satisfies the BSDE (\ref{BSDEgen}). Moreover, from the
strong convergence in $\mathbf{L}^{\mathbf{1}}(\tilde\mu)$ of $(U^n)$
to $U$, and
the Lipschitz condition on $h$, we have
\[
\mathbf{E}\int_{0}^T\int_{E} h^-(U^{n}_{s}(e),e)\lambda(de)\,ds
\rightarrow\mathbf{E}\int_{0}^T\int_{E}h^-(U_{s}(e),e)\lambda(de)\,ds
\]
as $n$ goes to infinity. Since $K_T^n = n \int_{0}^T\int_{E}
h^-(U^{n}_{s}(e),e)\lambda(de)\,ds$ is bounded in $\mathbf{L}^{\mathbf
{2}}(\bolds\Omega
,\bolds{\mathcal{F}}
_{\mathbf{T}})$, this implies
\[
\mathbf{E}\int_{0}^T\int_{E}h^-(U_{s}(e),e)\lambda(de)\,ds = 0
\]
and so the constraint (\ref{hcons}) is satisfied. Hence, $(Y,Z,K,U)$ is
a solution to the constrained BSDE (\ref{BSDEgen}) and (\ref{hcons}), and
by Lemma
\ref{lemcompbor}, $Y = \lim Y^n$ is the minimal solution. The
uniqueness of $Z$ follows by identifying the Brownian parts and the
finite variation parts, and then the uniqueness of $(U,K)$ is obtained
by identifying the predictable parts and by recalling that the jumps of
$\mu$ are inaccessible.

Finally, in the case $h(u,e) = -u$, the process
\[
\bar K_t = K_t - \int_0^t \int_E U_s(e) \mu(ds,de),\qquad 0 \leq
t\leq T,
\]
lies in $\mathbf{A}^{\mathbf{2}}$, and the triple $(Y,Z,\bar K)$ is
solution to (\ref
{BSDEneg}). Again, by Lemma \ref{lemcompbor},
this shows that $Y$ is the minimal solution to (\ref{BSDEgen}) and to
(\ref{BSDEneg}). The uniqueness of $(Y,Z,\bar K)$ is immediate by
identifying the Brownian part and the finite variation part.
\end{pf}
\begin{Remark} \label{remstrong}
From the estimate (\ref{estipeng}), it is clear that once
the process $K$ is continuous, that is, $\Delta K_t = 0$, then
$(Z^n,U^n)$ converges strongly to $(Z,U)$ in $\mathbf{L}^{\mathbf
{2}}(\mathbf{W})\times\mathbf{L}^{\mathbf{2}}(\tilde\mu)$. This occurs
in reflected
BSDEs as in \cite{elketal97} or \cite{hamouk03}; see also
Remark~\ref{remconvcont}. In the case of constraints on jump
component $U$ as in (\ref{BSDEgen}) and (\ref{hcons}), the situation
is more complicated, and the process $K$ is in general only
predictable. The same feature also occurs for constraints on $Z$
as in \cite{pen99}. To overcome this difficulty, we use the
estimations (\ref{Keps}) 
of the contribution of the
jumps of $K$, which allow us to obtain the strong convergence of
$(Z^n,U^n)$ in $\mathbf{L}^{\mathbf{p}}(\mathbf{W})\times\mathbf
{L}^{\mathbf{p}}(\tilde\mu)$ for $p \in[1,2)$. Finally, notice that
for the minimal solution
$(Y,Z,\tilde K)$ to the BSDE (\ref{BSDEneg}), the process $\tilde
K$ is not predictable.
\end{Remark}

\subsection{The case of impulse control}\label{sec33}

In the impulse control case [i.e., $f$ and $c$ depend only on $X$
and $h(u,e)=-u$], we have seen in Theorem \ref{theonegjump}
that the minimal solution to our constrained BSDE has the
following functional explicit representation:
\[
Y_t =
\mathop{\operatorname{ess}\sup}_{\nu\in\mathcal{V}} \mathbf{E}^\nu\biggl[ g(X_T) +
\int_t^T
f(X_s) \,ds + \int_t^T \int_E c(X_{s^-},e) \mu(ds,de) \Big|
\mathcal{F}_t \biggr].
\]
In this case, we also have a functional
explicit representation of the solution $Y^n$ to the penalized
BSDE (\ref{BSDEpen}),
%
%
\setcounter{equation}{21}
\renewcommand{\theequation}{\arabic{section}.\arabic{equation}}
\begin{eqnarray}\label{Ynimpulse}
Y^n_{t} &= & \mathop{\operatorname{ess}\sup}_{\nu\in\mathcal{V}_{n}}\mathbf
{E}^\nu\biggl[
g(X_{T})+\int_{t}^Tf(X_{s})\,ds\nonumber\\[-8pt]\\[-8pt]
& &\hspace*{46.5pt}{} + \int_{t}^T\int_{E}c(X_{s^-})\mu(ds,de)
\Big|\mathcal{F}_{t} \biggr],\nonumber
\end{eqnarray}
where $\mathcal{V}_{n}= \{
\nu\in\mathcal{V};\nu_{s}(e)\leq n\ \forall(s,e)\in[0,T]\times E$
a.s.$\}$. Indeed, denote by $\bar Y^n$ the right-hand side of
(\ref{Ynimpulse}). By writing that $(Y^n,Z^n,U^n)$ is the solution
of the penalized BSDE (\ref{BSDEpen}), taking the expectation
under $\mathbf{P}^\nu$, for $\nu\in\mathcal{V}_{n}$, and recalling
that $W$ is
a $\mathbf{P}^\nu$-Brownian motion, and $\nu\lambda(de)$ is the
intensity measure of $\mu$ under~$\mathbf{P}^\nu$, we obtain
%
%
\begin{eqnarray}\label{decompYimpulse}\qquad
Y^n_{t} & = & \mathbf{E}^\nu\biggl[ g(X_{T})+\int
_{t}^Tf(X_{s})\,ds+\int
_{t}^T\int_{E}c(X_{s^-},e)\mu(ds,de) \Big|\mathcal{F}_{t} \biggr]
\nonumber\\[-8pt]\\[-8pt]
& &{} + \mathbf{E}^\nu\biggl[\int_{t}^T\int_{E} \{
n[U_{s}^n(e)]_{+}-\nu_{s}(e)U_{s}^n(e) \}\lambda(de)\,ds
\Big|\mathcal{F}_{t} \biggr].\nonumber
\end{eqnarray}
Since this equality holds for any $\nu\in\mathcal{V}_n$, and
observing that
$n[U_{s}^n(e)]_{+}-\nu_{s}(e)U_{s}^n(e) \geq0$, for all $\nu\in
\mathcal{V}_{n}$, we have
%
%
\begin{equation}\label{decY2}\qquad
\bar Y_t^n \leq Y^n_{t} \leq \tilde{Y}^n_{t}
+ \mathbf{E}^\nu\biggl[\int_{t}^T\int_{E} \{
n[U_{s}^n(e)]_{+}-\nu
_{s}(e)U_{s}^n(e) \}\lambda(de)\,ds \Big|\mathcal{F}_{t}\biggr].
\end{equation}
Let us now consider the family $(\nu^\varepsilon)_{\varepsilon}$ of
$\mathcal{V}_{n}$ defined by
\[
\nu^\varepsilon_{s}(e) =
\cases{n,
&\quad if $U_{s}^n(e)>0$,\cr
\varepsilon, &\quad otherwise.}
\]
Then, by using the same argument as in the proof of Lemma \ref
{lemcompbor}, we show that
\[
\mathbf{E}^{\nu^\varepsilon} \biggl[\int_{t}^T\int_{E} \{
n[U_{s}^n(e)]_{+}-\nu_{s}(e)U_{s}^n(e) \}\lambda(de)\,ds
\Big|\mathcal{F}_{t} \biggr]
\rightarrow0 \qquad\mbox{as }\varepsilon\rightarrow0,
\]
which proves with (\ref{decY2}) that $Y^n_{t} = \bar{Y}^n_{t}$.

The representation (\ref{Ynimpulse}) has a nice interpretation. It
means that the value function of an impulse control problem can be
approximated by the value function of the same impulse control problem
but with strategies whose numbers of orders are bounded \textit{on
average} by $nT\lambda(E)$. This has to be compared with the classical
approximation by iterated optimal stopping problems, where the $n$th
iteration corresponds to the value of the same impulse control problem
but where the number of orders is smaller than $n$.
The numerical advantage of the penalized approximation is that it does
not require iterations.

\section{Relation with quasi-variational inequalities} \label{secrelQVI}

In this section, we show that minimal solutions to the
jump-constrained BSDEs provide a probabilistic representation of
solutions to parabolic QVIs of the form
%
%
\begin{eqnarray}\label{QVIgen}
\min\biggl[ - \frac{\partial v}{\partial t} - \mathcal{L}v -
f(\cdot,v,\sigma^{\intercal}D_xv) ,
\inf_{e\in
E} h( \mathcal{H}^e v - v,e) \biggr] = 0\nonumber\\[-8pt]\\[-8pt]
\eqntext{\mbox{on }
[0,T)\times\mathbb{R}^d,}
\end{eqnarray}
where $\mathcal{L}$ is the second-order local operator
\[
\mathcal{L}v(t,x)
= \langle b(x), D_x v(t,x)\rangle+ \tfrac{1}{2}\operatorname{
tr}(\sigma\sigma^{\intercal}(x)D_x^2v(t,x))
\]
and $\mathcal{H}^e$, $e \in E$, are the nonlocal operators
\[
\mathcal{H}^e v (t,x)= v\bigl(t,x+\gamma(x,e)\bigr) + c(x,v(t,x),\sigma
^{\intercal}(x)D_xv(t,x),e).
\]
For such nonlocal operators, we denote for $q \in\mathbb{R}^d$
\[
\mathcal{H}^e[t,x,q,v] = v\bigl(t,x+\gamma(x,e)\bigr) +
c(x,v(t,x),\sigma^{\intercal}(x)q,e).
\]

Note that when $h(u)$ does not depend on $e$, and since it is
nonincreasing in $u$, the QVI (\ref{QVIgen}) may be written
equivalently in
\[
\min\biggl[ - \frac{\partial v}{\partial t} -
\mathcal{L}v -
f(\cdot,v,\sigma^{\intercal}D_xv) , h( \mathcal{H}v - v)\biggr] = 0
\qquad\mbox{on } [0,T)\times\mathbb{R}^d,
\]
with $\mathcal{H}v = \sup_{e\in E}\mathcal{H}^ev$. In particular, this
includes the case of
QVI associated to impulse controls for $h(u) = -u$, and $f$,
$c$ independent of $y,z$.

We shall use the penalized parabolic integral partial differential
equation (IPDE) associated to the penalized BSDE (\ref{BSDEpen}), for
each $n \in\mathbb{N}$,
%
%
\begin{eqnarray}\label{IPDE}
&& - \frac{\partial v_n}{\partial t} - \mathcal{L}v_n -
f(\cdot,v_n,\sigma^{\intercal}D_x v_n) \nonumber\\[-8pt]\\[-8pt]
&&\qquad{} - n \int_E
h^-
( \mathcal{H}^ev_n -
v_n ,e ) \lambda(de) = 0\qquad \mbox{on } [0,T)\times
\mathbb{R}^d.\nonumber
\end{eqnarray}

To complete the PDE characterization of the function $v$, we need to
provide a suitable boundary condition. In general, we cannot expect
to have $v(T^- , \cdot) = g$, and we shall consider the relaxed boundary
condition given by the equation
%
%
\begin{equation}
\label{bound-cond}\qquad
\min\Bigl[ v(T^- , \cdot) - g , \inf_{e\in E} h\bigl(
\mathcal{H}^e v(T^- , \cdot)
- v(T^- , \cdot),e\bigr) \Bigr] = 0\qquad \mbox{on } \mathbb{R}^d.
\end{equation}

In the sequel, we shall assume in addition to the conditions of
Section \ref{secgeneral} that the functions $\gamma$, $f$, $c$
and $h$ are continuous with respect to all their arguments.

\subsection{Viscosity properties}\label{sec41}

Solutions of (\ref{QVIgen}), (\ref{IPDE}) and (\ref{bound-cond}) are
considered in the (discontinuous) viscosity sense, and it will be
convenient in the sequel to define the notion of viscosity solutions
in terms of sub- and super-jets. We refer to \cite{ish93,tanyon93}
and more recently to the book \cite{okssul07}
for the notion of viscosity solutions to QVIs.
For a locally bounded function $u$
on $[0,T]\times\mathbb{R}^d$, we define its lower semicontinuous (lsc in
short) $u_*$, and upper semicontinuous (usc in short) envelope $u^*$
by
\[
u_*(t,x) = \liminf_{\tiny{(t',x')\rightarrow(t,x),
t'<T}} u(t',x'),\qquad u^*(t,x) =
\limsup_{\tiny{(t',x')\rightarrow(t,x), t'<T}} u(t',x').
\]
\begin{Definition}[(Subjets and superjets)]
(i) For a function $u \dvtx[0,T]\times\mathbb{R}^d
\rightarrow
\mathbb{R}$, lsc (resp., usc), we denote by $J^-u(t,x)$ the parabolic
subjet [resp., $J^+u(t,x)$ the parabolic superjet] of $u$ at
$(t,x) \in[0,T]\times\mathbb{R}^d$, as the set of triples $(p,q,M) \in
\mathbb{R}\times\mathbb{R}^d\times\mathbb{S}^d$ satisfying
\begin{eqnarray*}
u(t',x') & \geq&(\mbox{resp.,} \leq)\  u(t,x) + p(t'-t) + \langle
q,x'-x\rangle+
\tfrac{1}{2}\langle x'-x, M(x'-x)\rangle\\
& &\hspace*{42.7pt}{} + o(|t'-t| + |x'-x|^2).
\end{eqnarray*}

(ii) For a function $u \dvtx[0,T)\times\mathbb{R}^d \rightarrow
\mathbb{R}$,
lsc (resp., usc), we denote by $\bar J^-u(t,x)$ the parabolic limiting
subjet [resp.,
$\bar J^+u(t,x)$ the parabolic limiting superjet] of $u$ at $(t,x) \in
[0,T]\times\mathbb{R}^d$,
as the set of triples $(p,q,M) \in\mathbb{R}\times\mathbb
{R}^d\times\mathbb{S}^d$
such that
\begin{eqnarray}
(p,q,M) = \lim_{n} (p_n,q_n,M_n), \qquad (t,x) = \lim_n
(t_n,x_n) \nonumber\\
\eqntext{\mbox{with } (p_n,q_n,M_n) \in J^-u(t_n,x_n) \mbox{ [resp.,
$J^+u(t_n,x_n)$]},}\\
\eqntext{\displaystyle
u(t,x) = \lim_n u(t_n,x_n).}
\end{eqnarray}
\end{Definition}

We now give the definition of viscosity solutions to (\ref{QVIgen}),
(\ref{IPDE}) and (\ref{bound-cond}).
\begin{Definition}[{[Viscosity solutions to (\ref{QVIgen})]}]
(i) A function $u$, lsc (resp., usc) on $[0,T)\times\mathbb{R}^d$,
is called a viscosity supersolution (resp., subsolution) to
(\ref{QVIgen}) if for each $(t,x) \in[0,T)\times\mathbb{R}^d$,
and any
$(p,q,M) \in\bar J^- u(t,x)$ [resp., $\bar J^+u(t,x)$], we have
\begin{eqnarray*}
&&\min\biggl[ - p - \langle b(x),q \rangle- \frac{1}{2}\operatorname{
tr}(\sigma\sigma
^{\intercal}(x)M) \\
&&\hspace*{20.7pt}{} - f(x,u(t,x),\sigma^{\intercal}(x)q) ,
\inf_{e\in E} h\bigl( \mathcal{H}^e [t,x,q,u] - u(t,x),e\bigr) \biggr]  \geq
\mbox{(resp.,  $\leq$) }  0.
\end{eqnarray*}

(ii) A locally bounded function on $[0,T)\times\mathbb{R}^d$ is
called a
viscosity solution to (\ref{QVIgen}) if $u_*$ is a viscosity supersolution
and $u^*$ is a viscosity subsolution to (\ref{QVIgen}).
\end{Definition}
\begin{Definition}[{[Viscosity solutions to (\ref{IPDE})]}]
(i) A function $u$, lsc (resp., usc) on $[0,T)\times\mathbb{R}^d$,
is called a viscosity supersolution (resp., subsolution) to (\ref{IPDE})
if for
each $(t,x) \in[0,T)\times\mathbb{R}^d$, and any $(p,q,M) \in\bar
J^- u(t,x)$ [resp., $\bar J^+u(t,x)$], we have
\begin{eqnarray*}
&&- p - \langle b(x),q\rangle- \frac{1}{2}\operatorname{tr}(\sigma\sigma
^{\intercal}(x)M)
- f(x,u(t,x),\sigma^{\intercal}(x)q)\\
&&\qquad{}- n \int_E h^-\bigl( \mathcal{H}^e [t,x,q,u] - u(t,x),e\bigr) \lambda(de)
\geq
\mbox{(resp.,  $\leq$) } 0.
\end{eqnarray*}

(ii) A locally bounded function $u$ on $[0,T)\times\mathbb{R}^d$ is
called a
viscosity solution to (\ref{IPDE}) if $u_*$ is a viscosity
supersolution and $u^*$ is a viscosity subsolution to (\ref{IPDE}).
\end{Definition}
\begin{Definition}[{[Viscosity solutions to (\ref{bound-cond})]}]
(i) A function $u$, lsc (resp., usc) on $[0,T]\times\mathbb{R}^d$,
is called a viscosity supersolution (resp., subsolution) to (\ref
{bound-cond}) if for
each $x \in\mathbb{R}^d$, and any $(p,q,M) \in\bar J^- u(T,x)$
[resp., $\bar J^+u(T,x)$], we have
\[
\min\Bigl[ u(T,x)-g(x) ,
\inf_{e\in E} h\bigl( \mathcal{H}^e [T,x,q,u] - u(T,x),e\bigr) \Bigr] \geq
\mbox{(resp., $\leq$) } 0.
\]

(ii) A locally bounded function $u$ on $[0,T]\times\mathbb{R}^d$ is
called a
viscosity solution to (\ref{bound-cond}) if $u_*$ is a viscosity supersolution
and $u^*$ is a viscosity subsolution to (\ref{bound-cond}).
\end{Definition}
\begin{Remark}
An equivalent definition of viscosity super and subsolution to~(\ref
{bound-cond}), which shall be used later, is the following in terms of
test functions:
a~function $u$, lsc (resp., usc) on $[0,T]\times\mathbb{R}^d$, is
called a
viscosity supersolution (resp., subsolution) to (\ref{bound-cond}) if for
each $(t,x) \in[0,T)\times\mathbb{R}^d$, and any $\varphi\in
C^{1,2}([0,T]\times\mathbb{R}^d)$ such that $(t,x)$ is a
minimum (resp., maximum) global of $u-\varphi$, we have
\[
\min\Bigl[ u(T,x)-g(x) ,
\inf_{e\in E} h\bigl( \mathcal{H}^e [T,x,D_x\varphi(T,x),u] -
u(T,x),e\bigr) \Bigr]
\geq(\mbox{resp.}, \leq)\ 0.
\]
We have similar equivalent definitions of viscosity super and
subsolution to (\ref{QVIgen}) in terms of test functions.
\end{Remark}

We slightly strengthen assumption (\ref{assuH1}) or (\ref{assuH2}) by
\renewcommand{\theequation}{H1$'$}
\begin{eqnarray}\label{assuH1prime}
\\
\nonumber\\[-36.5pt]
&\begin{tabular}{p{325pt}}
There exists a quadruple $(\tilde
Y,\tilde Z,\tilde K) \in\bolds{\mathcal{S}}^{\mathbf{2}}\times\mathbf
{L}^{\mathbf{2}}(\mathbf{W})\times\mathbf{A}^{\mathbf{2}}$ satisfying
(\ref{BSDEneg}), with $\tilde Y_t = \tilde v(t,X_t)$,
$0\leq t\leq T$, for some function deterministic $\tilde v$ satisfying
a linear growth condition
\end{tabular}\hspace*{-36pt}&
\nonumber\\
&\displaystyle
\sup_{(t,x)\in[0,T]\times\mathbb{R}^d}\frac{|\tilde
v(t,x)|}{1+|x|}<+\infty.&\nonumber\\[-30pt]\nonumber
\end{eqnarray}
\renewcommand{\theequation}{H2$'$}
\begin{eqnarray}\label{assuH2prime}
\\
\nonumber\\[-36.5pt]
&\begin{tabular}{p{325pt}}
There exists a quadruple $(\tilde
Y,\tilde Z,\tilde K,\tilde U) \in\bolds{\mathcal{S}}^{\mathbf{2}}\times
\mathbf{L}^{\mathbf{2}}(\mathbf{W})\times\mathbf{L}^{\mathbf{2}}(\tilde
\mu)\times
\mathbf{A}^{\mathbf{2}}$ satisfying (\ref{BSDEgen}) and (\ref{hcons}), with
$\tilde Y_t = \tilde v(t,X_t)$,
$0\leq t\leq T$, for some function deterministic $\tilde v$ satisfying
a linear growth condition
\end{tabular}\hspace*{-36pt}&
\nonumber\\
&\displaystyle
\sup_{(t,x)\in[0,T]\times\mathbb{R}^d}\frac{|\tilde
v(t,x)|}{1+|x|}<+\infty.&\nonumber
\end{eqnarray}

Under assumption (\ref{assuH1prime}) [resp., (\ref{assuH2prime})], there esists for each
$(t,x) \in[0,T]\times\mathbb{R}^d$ a unique minimal solution $\{
(Y_{s}^{t,x},Z_{s}^{t,x},U_{s}^{t,x},K_{s}^{t,x}),t\leq s\leq T\}$ to
(\ref{BSDEgen}) and (\ref{hcons}) [resp., (\ref{BSDEneg}) and (\ref{negjump})]
with $X = \{X_{s}^{t,x},t\leq s\leq T\}$, the solution to (\ref{defXjump})
starting from $x$ at time $t$.
We can then define the (deterministic) function $v\dvtx[0,T]\times
\mathbb{R}
^d\rightarrow\mathbb{R}$ by
%
%
\setcounter{equation}{3}
\renewcommand{\theequation}{\arabic{section}.\arabic{equation}}
\begin{equation}\label{defvY}
v(t,x) := Y_t^{t,x},\qquad (t,x) \in[0,T]\times\mathbb{R}^d.
\end{equation}
%
Similarly, we define the function
%
%
\begin{equation}\label{defvnY}
v_n(t,x) := Y_t^{n,t,x}, \qquad (t,x) \in[0,T]\times\mathbb{R}^d,
\end{equation}
where $\{(Y_s^{n,t,x},Z_s^{n,t,x},U_s^{n,t,x}(\cdot)), t\leq s\leq T\}$ is
the unique solution to (\ref{BSDEpen}) with $X_s = X_s^{t,x}$,
$t\leq s\leq T$.

We first have the following identification.
\begin{Proposition}
The function $v$ links the processes $Y^{t,x}$ and $X^{t,x}$ by the relation
%
%
\begin{equation}\label{idv}
Y^{t,x}_{\theta} = v (\theta,X_{\theta}^{t,x} )\qquad
\mbox{for all stopping time } \theta\mbox{ valued in } [t,T].
\end{equation}
\end{Proposition}
\begin{pf}
From the Markov property of the jump-diffusion process $X$, and
uniqueness of a solution $Y^n$ to the BSDE (\ref{BSDEpen}), we have
(see, e.g., \cite{barbucpar97})
%
%
\begin{equation}\label{idvn}
Y^{t,x,n}_{s} = v_{n} (s,X_{s}^{t,x} ),\qquad t\leq
s\leq T.
\end{equation}
From Section \ref{sec3}, we know that $v$ is the pointwise limit of $v_{n}$.
Moreover, by (\ref{YnL2}),
$Y^{t,x,n}_{\theta}$ converges to $Y^{t,x}_{\theta}$ as $n$ goes to
infinity, for all stopping time $\theta$ valued in $[t,T]$. We then
obtain the required relation by passing to the limit in (\ref{idvn}).
\end{pf}
\begin{Remark} \label{remvgrowth}
Assumption (\ref{assuH2prime}) [or (\ref{assuH1prime}) which is weaker than
(\ref{assuH2prime}) in the case $h(u,e) = -u$]
ensures that the function $v$ in (\ref{defvY}) satisfies a linear
growth condition, and is in particular locally bounded.
Indeed, from (\ref{estiYngrowth}) and by passing to the limit by
Fatou's lemma for $v(t,x) = Y_t^{t,x} = \lim Y_t^{n,t,x}$, we have
\begin{eqnarray*}
{\sup_{t\in[0,T]}} |v(t,x)|^2 & \leq& C \biggl( 1 + \mathbf{E}|g(X_T^{t,x})|^2
+ \mathbf{E} \biggl[\int_t^T |X_s^{t,x}|^2 \,dt \biggr]\\
& &\hspace*{65.5pt}{} + \mathbf{E} \Bigl[{\sup_{s\in[t,T]}}|\tilde
v(s,X_s^{t,x})|^2 \Bigr] \biggr).
\end{eqnarray*}
The result follows from the standard estimate
\[
\mathbf{E} \Bigl[{\sup_{t\leq s \leq T}} |X_s^{t,x}|^2 \Bigr] \leq C(1+|x|^2)
\]
and the linear growth conditions on $g$ and $\tilde v$.
\end{Remark}

The relation between the penalized BSDE (\ref{BSDEpen}) and the
penalized IPDE (\ref{IPDE}) is well known from the results of \cite
{barbucpar97}.
Although our framework does not fit exactly into the one of \cite
{barbucpar97}, by mimicking closely the arguments in this paper and
using comparison theorem in \cite{roy06}, we obtain the following result.
\begin{Proposition} \label{provn}
The function $v_n$ in (\ref{defvnY}) is a continuous viscosity solution
to (\ref{BSDEpen}).
\end{Proposition}

By adapting stability arguments for viscosity solutions to our context,
we now prove the viscosity property of the function $v$ to (\ref{QVIgen}).
We shall assume that the support of $\lambda$ is the whole space $E$, that is,
\renewcommand{\theequation}{HE}
\begin{eqnarray}\label{assuHE}
&&\mbox{$\forall e \in E\qquad \exists\mathcal{O}$ open
neighborhood of $e$, s.t. $\lambda(\mathcal{O}) > 0$.}
\end{eqnarray}
\begin{Theorem}
Under (\ref{assuH2prime}) [or (\ref{assuH1prime}) in the case: $h(u,e) = -u$], and
(\ref{assuHE}), the function $v$ in (\ref{defvY}) is a viscosity solution
to (\ref{QVIgen}).
\end{Theorem}
\begin{pf}
From the results of the previous section, we know that $v$ is the
pointwise limit
of the nondecreasing sequence of functions $(v_n)$. By continuity of
$v_n$, we then have (see, e.g., \cite{bar94}, page 91):
%
%
\setcounter{equation}{7}
\renewcommand{\theequation}{\arabic{section}.\arabic{equation}}
\begin{eqnarray}\label{vinf}
v &=& v_* = \lim_{n\rightarrow\infty}{
\inf}_* v_n\hspace*{70pt} \nonumber\\[-8pt]\\[-8pt]
\eqntext{\mbox{where }
\displaystyle
\lim_{n\rightarrow\infty}{\inf}_* v_n(t,x) := \mathop{\liminf_{
n\rightarrow\infty}}_{t'\rightarrow t,x'\rightarrow x
} v_n(t',x'),}
\\
\label{vsup}
v^* &=& \lim_{n\rightarrow\infty}{\sup}^*v_n\hspace*{70pt} \nonumber\\[-8pt]\\[-8pt]
\eqntext{\mbox{where }
\displaystyle
\lim_{n\rightarrow\infty}{\sup}^* v_n(t,x) := \mathop{\limsup
_{n\rightarrow\infty}}_{t'\rightarrow t,x'\rightarrow x} v_n(t',x').}
\end{eqnarray}

(i) We first show the viscosity supersolution property for $v = v_*$.
Let $(t,x)$ be a point in
$[0,T)\times\mathbb{R}^d$, and $(p,q,M) \in\bar J^-v(t,x)$. By
(\ref{vinf}) and Lemma 6.1 in \cite{craishlio92}, there exists sequences
\[
n_j \rightarrow\infty,\qquad (p_j,q_j,M_j) \in J^-
v_{n_j}(t_j,x_j),
\]
such that
%
%
\begin{equation}\label{pjqj2}
(t_j,x_j,v_{n_j}(t_j,x_j),p_j,q_j,M_j) \rightarrow(t,x,v(t,x),p,q,M).
\end{equation}
We also have by definition of $v = v_*$ and continuity of $\gamma$:
%
%
\begin{equation}\label{limvgam}
v\bigl(t,x+\gamma(x,e)\bigr) \leq\liminf_{j\rightarrow\infty}
v_{n_j}\bigl(t_j,x_j+\gamma(x_j,e)\bigr)\qquad \forall e \in E.
\end{equation}
%
Moreover, from the viscosity supersolution property for $v_{n_j}$, we
have for all $j$
%
%
\begin{eqnarray} \label{supervnj}\qquad
&&- p_j - \langle b(x_j),q_j \rangle- \frac{1}{2}\operatorname{tr}(\sigma
\sigma
^{\intercal}(x_j)M_j) - f(x_j,v_{n_j}(t_j,x_j),\sigma^{\intercal
}(x_j)q_j)
\nonumber\\[-8pt]\\[-8pt]
&&\qquad{} - n_j \int_E h^-\bigl( \mathcal{H}^e [t_j,x_j,q_j,v_{n_j}] - v_{n_j}(t_j,x_j),e\bigr)
\lambda(de) \geq 0.\nonumber
\end{eqnarray}
Let us check that the following inequality holds:
%
%
\begin{equation}
\label{intersur}
\inf_{e\in E} h\bigl( \mathcal{H}^e[t,x,q,v] - v(t,x),e\bigr)
\geq0.
\end{equation}
We argue by contradiction, and assume there exists
some $e_0 \in E$ s.t.
\[
h\bigl(v\bigl(t,x+\gamma(x,e_0)\bigr) +
c(x,v(t,x),\sigma^{\intercal}(x) q,e_0)- v(t,x),e_0\bigr) < 0.
\]
Then, by continuity of $\sigma$, $h$, $\gamma$, $c$ in all their
variables, (\ref{pjqj2}), (\ref{limvgam}) and the nonincreasing
property of $h$, one may find some $\varepsilon> 0$ and some open
neighborhood $\mathcal{O}_0$ of $e_0$ such that for all $j$ large
enough
\[
h\bigl(v_{n_j}\bigl(t_j,x_j+\gamma(x_j,e)\bigr) +
c(x_j,v_{n_j}(t_j,x_j),\sigma^{\intercal}(x_j) q_j,e)-
v_{n_j}(t_j,x_j),e\bigr) \leq- \varepsilon 
\]
for all $e \in\mathcal{O}_0$. Since the support of $\lambda$ is $E$,
this implies
\[
\int_E h^-\bigl(\mathcal{H}^e(t_j,x_j,q_j,v_{n_j}) - v_{n_j}(t_j,x_j),e\bigr)
\lambda(de) \geq\varepsilon\lambda(\mathcal{O}_0) > 0.
\]
By
sending $j$ to infinity into (\ref{supervnj}), we get the required
contradiction. On the other hand, by (\ref{supervnj}), we have
\[
- p_j - \langle b(x_j),q_j \rangle- \tfrac
{1}{2}\operatorname{tr}(\sigma\sigma^{\intercal}(x_j)M_j) -
f(x_j,v_{n_j}(t_j,x_j),\sigma^{\intercal}(x_j)q_j) \geq0,
\]
so
that by sending $j$ to infinity,
\[
- p -\langle b(x),q\rangle- \tfrac{1}{2}\operatorname{tr}
(\sigma\sigma^{\intercal}(x)M) - f(x,v(t,x),\sigma^{\intercal
}(x)q) \geq0,
\]
which proves, together with (\ref{intersur}), that $v$ is a
viscosity supersolution to (\ref{QVIgen}).

(ii) We conclude by showing the viscosity subsolution
property for $v^*$. Let $(t,x)$ a point in $[0,T)\times\mathbb{R}^d$, and
$(p,q,M) \in\bar J^+v^*(t,x)$ such that
%
%
\begin{equation}\label{hce}
\inf_{e\in E} h\bigl( \mathcal{H}^e[t,x,q,v^*] - v^*(t,x),e\bigr) > 0.
\end{equation}
From (\ref{vsup}) and Lemma 6.1 in \cite{craishlio92}, there
exist sequences
\[
n_j \rightarrow
\infty,\qquad (p_j,q_j,M_j) \in J^+ v_{n_j}(t_j,x_j),
\]
such that
%
%
\begin{equation}\label{pjqj}
(t_j,x_j,v_{n_j}(t_j,x_j),p_j,q_j,M_j) \rightarrow
(t,x,v^*(t,x),p,q,M).
\end{equation}
By continuity of the functions
$c,\gamma$ and the definition of $v^*$, we also have
%
%
\begin{equation}
\label{vnjgam} \limsup_{j\rightarrow\infty}
\mathcal{H}^e[t_j,x_j,q_{j},v_{n_j}] \leq\mathcal
{H}^e[t,x,q,v^*]\qquad\forall e \in E.
\end{equation}
Now, from the viscosity subsolution
property for $v_{n_j}$, we have for all $j$
%
%
\begin{eqnarray} \label{sousolvnj}\qquad
&&- p_j - \langle
b(x_j),q_j\rangle- \frac{1}{2}\operatorname{tr}(\sigma\sigma^{\intercal
}(x_j)M_j)
-f(x_j,v_{n_j}(t_j,x_j),\sigma^{\intercal}(x_j)q_j)\nonumber\\[-8pt]\\[-8pt]
&&\qquad{} - n_j \int_E h^-\bigl( \mathcal{H}^e [t_j,x_j,q_j,v_{n_j}] - v_{n_j}(t_j,x_j),e\bigr)
\lambda(de) \leq 0.\nonumber
\end{eqnarray}
From (\ref{hce}) (which is uniform in $e\in E$), (\ref{pjqj}) and (\ref
{vnjgam}), continuity
assumptions on $h,c$ and the nonincreasing property of $h$, we
have for $j$ large enough
\[
h\bigl( \mathcal{H}^e [t_j,x_j,q_j,v_{n_j}] -
v_{n_j}(t_j,x_j),e\bigr) > 0\qquad \forall e \in E,
\]
and so
\[
\int_E h^-\bigl( \mathcal{H}^e [t_j,x_j,q_j,v_{n_j}] -
v_{n_j}(t_j,x_j),e\bigr) \lambda(de) = 0.
\]
Hence, by taking the
limit as $j$ goes to infinity, into (\ref{sousolvnj}), we conclude
that
\[
- p - \langle b(x),q \rangle- \tfrac{1}{2}\operatorname{tr}(\sigma\sigma
^{\intercal}(x)M) -
f(x,v^*(t,x),\sigma^{\intercal}(x)q) \leq0,
\]
which shows the
viscosity subsolution property for $v^*$ to (\ref{QVIgen}).
\end{pf}

We next turn to the boundary condition.
\begin{Theorem} \label{probound}
Under (\ref{assuH2prime}) [or (\ref{assuH1prime}) in the case: $h(u,e) = -u$] and
(\ref{assuHE}), the function $v$ in (\ref{defvY}) is a viscosity solution
to (\ref{bound-cond}).
\end{Theorem}

In order to deal with the possible jump at the terminal condition, we
need the following dynamic programming characterization of the minimal solution.
\begin{Lemma} \label{lemprog}
Let $(t,x) \in[0,T)\times\mathbb{R}^d$, and
$(Y^{t,x},Z^{t,x},U^{t,x},K^{t,x})$ be a minimal solution to (\ref
{BSDEgen}) and (\ref{hcons}) on $[t,T]$ with $X_s = X_s^{t,x}$.
Then for any stopping time $\theta$ valued in $[t,T]$,
$(Y^{t,x}_{s},Z^{t,x}_{s},U^{t,x}_{s},K^{t,x}_{s})_{s\in[t,\theta]}$
is a minimal solution to
%
%
\setcounter{equation}{17}
\renewcommand{\theequation}{\arabic{section}.\arabic{equation}}
\begin{eqnarray} \label{EDSRPPD}
Y_{s} & = & v(\theta,X_{\theta}^{t,x})+\int_{s}^\theta
f(X^{t,x}_{r},Y_{r},Z_{r})\,dr\nonumber\\
& &{} + K_\theta^{t,x} - K_s^{t,x} 
-\int_{s}^\theta\langle Z_{r},dW_{r}\rangle
\\
& &{} -\int_{s}^\theta\int_{E} \bigl(
U_{r}(e)-c(X_{r^-}^{t,x},Y_{r^-},Z_{r},e) \bigr) \mu(dr,de) \nonumber
\end{eqnarray}
with
%
%
\begin{equation}\label{hconsPPD}
h(U_{s}(e),e)\geq0,\qquad d\mathbf{P}\otimes dt\otimes
\lambda
(de),\mbox{ a.e. on } \Omega\times[t,\theta] \times E.
\end{equation}
\end{Lemma}
\begin{pf}
Notice first from (\ref{idv}) that
$(Y_{s}^{t,x},Z_{s}^{t,x},U_{s}^{t,x},K_{s}^{t,x})_{s\in[t,\theta]}$
is solution to (\ref{EDSRPPD}) and (\ref{hconsPPD}).
Let $Y^1$ be the minimal solution on $[t,\theta]$ of (\ref
{EDSRPPD}) and (\ref{hconsPPD})
(the existence of a minimal solution in the case of a random terminal
time is obtained by similar arguments to those used in the case of a
deterministic terminal time). For each\vspace*{-1pt} $\omega\in\Omega$, there
exists a minimal solution $Y^{2,\omega}$ on $[\theta(\omega),T]$ to
(\ref{BSDEgen}) and (\ref{hcons}) with\vspace*{2pt} $X = \{X_{s}^{\theta(\omega
),X^{t,x}_{\theta(\omega)}(\omega)}, \theta(\omega)\leq s\leq T\}$.
We then have from the definition of $v$
that $Y^{2,\omega}_{\theta(\omega)}=v (\theta(\omega
),X^{t,x}_{\theta(\omega)}(\omega) )$ for
all $\omega\in\Omega$. By a measurable selection result (see, e.g.,
Theorem~82 in the Appendix to Chapter III in \cite{delmey75}),
there exists $Y^2 \in\-\bolds{\mathcal{S}}^{\mathbf{2}}$ such that
$\mathbf
{P}(d\omega)$
a.s., we have
$Y^{2}_{\theta(\omega)}(\omega)=Y^{2,\omega}_{\theta(\omega
)}=v (\theta(\omega),X^{t,x}_{\theta(\omega)}(\omega)
)$ and
$Y^{2}_{s}(\omega)=Y^{2,\omega}_{s}(\omega)$ for $s\in[\theta
(\omega),T]$.
We then define the process $\tilde{Y}$ by $\tilde{Y}|_{[t,\theta
]}=Y^1$ and $\tilde{Y}|_{(\theta,T]}=Y^2$.
Hence, $\tilde{Y}$ is a solution on $[t,T]$ to (\ref{BSDEgen}) and (\ref
{hcons}), which implies $\tilde{Y}\geq Y^{t,x}$. Moreover, since
$Y^{t,x}_{\theta} = v(\theta,X^{t,x}_{\theta})$, it follows that
$(Y^{t,x}_{s},Z^{t,x}_{s},U^{t,x}_{s},K^{t,x}_{s})_{s\in[t,\theta]}$
is a solution on $[t,\theta]$ to (\ref{EDSRPPD}) and (\ref{hconsPPD}).
Hence, $Y^1\leq Y^{t,x}$ on $[t,\theta]$, and therefore $Y^1= Y^{t,x}$
on $[t,\theta]$.
\end{pf}
\begin{pf*}{Proof of Theorem \protect\ref{probound}}
(i) We first prove the supersolution property of $v_*$ to
(\ref{bound-cond}).
Let $x\in\mathbb{R}^d$, and $(p,q,M) \in\bar J^-v_*(T,x)$. By the
same arguments as in (\ref{intersur}), we have
%
%
\setcounter{equation}{19}
\renewcommand{\theequation}{\arabic{section}.\arabic{equation}}
\begin{equation}\label{intersurT}
\inf_{e\in E}h\bigl(\mathcal{H}^{e}[T,x, q,v_*]-v_*(T,x),e\bigr) \geq0.
\end{equation}
Moreover, since the sequence of continuous functions $(v_{n})_{n}$ is
nondecreasing and $v_{n}(T,\cdot)=g$, we deduce that $v_*(T,\cdot) \geq g$,
which combined with (\ref{intersurT}), proves the viscosity
supersolution property for $v_*$ to (\ref{bound-cond}).

(ii) We next prove the subsolution property of $v^*$ to (\ref
{bound-cond}).
We argue by contradiction and assume that there exist $x_{0}\in\mathbb
{R}^n$, $\varphi\in C^{1,2}([0,T]\times\mathbb{R}^n)$
such that
%
%
\begin{equation}\label{varphimax}
0 = (v^*-\varphi)(T,x_{0}) = \max_{[0,T]\times\mathbb
{R}^d}(v^*-\varphi)
\end{equation}
and
\begin{eqnarray*}
&&\min\Bigl[ \varphi(T,x_{0})-g(x_{0}) ,
\inf_{e\in E}h\bigl(\mathcal{H}^{e}[T,x_{0},D_x\varphi(T,x_{0}),v^*]
-\varphi(T,x_{0}),e\bigr) \Bigr]\\
&&\qquad =: 2\varepsilon> 0.
\end{eqnarray*}
By the upper semicontinuity of $v^*$, the continuity of $\varphi$ and
its derivative, and the nonincreasing property of $h$, there exists an
open neighbohood $\mathcal{O}$ of $(T,x_{0})$ in $[0,T]\times\mathbb{R}^d$,
and $A,r>0$ such that for all
$(t,x,\alpha,\beta)\in\mathcal{O}\times(-A,A)\times B(0,r)$, we have
%
%
\begin{eqnarray}\label{cond0}
\varepsilon& \leq& \min\Bigl[\varphi(t,x)-\alpha
-g(x) ,
\nonumber\\
& &\hspace*{22.3pt} \inf_{e\in E}h \bigl(v^*\bigl(t,x+\gamma(x,e)\bigr)
\nonumber\\[-8pt]\\[-8pt]
& &\hspace*{47pt}{} +
c\bigl(x,\varphi(t,x)-\alpha,\sigma^{\intercal}(x)[D_x\varphi
(t,x)+\beta]\bigr) \nonumber\\
&&\hspace*{157.6pt}{} -
[\varphi(t,x)-\alpha],e \bigr) \Bigr]. \nonumber
\end{eqnarray}
Let $(t_{k},x_{k})_{k}$ be a sequence in $[0,T)\times\mathbb{R}^d$
such that
%
%
\begin{equation}
\label{convbord}
(t_{k},x_{k})\rightarrow(T,x_{0}) \quad\mbox{and}\quad
v(t_{k},x_{k})\rightarrow v^*(T,x_{0}).
\end{equation}
Fix then $\delta>0$ such that for $k$ large enough: $[t_{k},T]\times
B(x_{k},\delta)\subset\mathcal{O}$, and let us define the functions
$\varphi_{k}$ by
\[
\varphi_{k}(t,x) =
\varphi(t,x)+\zeta\frac{| x-x_{k}
|^2}{\delta^2}+C_{k}\phi\biggl(\frac{x-x_{k}}{\delta} \biggr)+\sqrt{T-t},
\]
where $0<\zeta<A\wedge\delta r$, $\phi\in C^2(\mathbb{R}^d)$
satisfies $\phi|_{\bar{B}(0,1)}\equiv0,\phi|_{\bar{B}(0,1)^{c}}>
0$ and $\lim_{|x|\rightarrow\infty}\frac{\phi(x)}{1+|x|}=\infty$,
and\vspace*{1pt} $C_{k}>0$ is a constant to be chosen below. By (\ref{varphimax}),
we notice that
\[
(v^*-\varphi_{k})(t,x) \leq-\zeta\qquad\mbox{for } (t,x)\in
[t_{k},T]\times\partial B(x_{k},\delta)
\]
and from the conditions on $\phi$, we can choose $C_{k}$ (large
enough) so that
%
%
\begin{equation}\label{majhdom}
(v^*-\varphi_{k})(t,x) \leq-\frac{\zeta}{2} \qquad\mbox{for
} (t,x)\in[t_{k},T]\times B(x_{k},\delta)^c.
\end{equation}
Since $\frac{\partial}{\partial
t}(\sqrt{T-t})\rightarrow-\infty$ as $t\nearrow T$, we have for
$k$ large enough
%
%
\begin{eqnarray} \label{cond1}
&& - \frac{\partial\varphi_k}{\partial t} - \mathcal{L}\varphi
_{k}(t,x)-f\bigl(x,\varphi
_{k}(t,x)-\alpha,\sigma^{\intercal}(x)D_x\varphi_{k}(t,x)\bigr)
\geq 0 \nonumber\\[-8pt]\\[-8pt]
\eqntext{\mbox{for } (t,x,\alpha)\in[t_{k},T)\times
B(x_{k},\delta)\times(-A+\zeta,A).}
\end{eqnarray}
Fix now $\alpha^*\in(0,A\wedge\frac{\zeta}{2}\wedge\varepsilon
)$, and let us denote
$\tau_{k} = \inf\{ s\geq t_{k};X_{s}^k\neq X_{s^-}^k
\}$,
$\theta_{k} = \inf\{ s\geq t_{k};X_{s}^k \notin
B(x_k,\delta) \}\wedge\tau_{k}\wedge T$ where $X^k=X^{t_{k},x_{k}}$.
Let us then define the quadruples $(Y^k, Z^k,U^k,K^k)$ on
$[t_{k},\theta_{k}]$ by
\begin{eqnarray*}
Y_{s}^k & = & [\varphi_{k}(s,X^k_{s})-\alpha^* ]\mathbf
{1}_{\{s\in[t_{k},\theta_{k})\}} +
v(\theta_{k},X^k_{\theta_{k}})\mathbf{1}_{\{s=\theta_{k}\}} ,\\
Z^k_{s} & = & \sigma^{\intercal}(X_{s^-}^k)D_x\varphi
_{k}(s,X^k_{s^-}) ,\\
U^k_{s}(e) & = &
v^*\bigl(s,X^k_{s^-}+\gamma(X^k_{s^-},e)\bigr) \\
& &{} +c\bigl(X^k_{s^-},\varphi_{k}(s,X^k_{s^-})-\alpha^*,
\sigma^{\intercal}(X^k_{s^-}) D_x\varphi_{k}(s,X^k_{s^-})\bigr) \\
& &{} - [\varphi_{k}(s,X^k_{s^-})-\alpha^*]
\end{eqnarray*}
and
\begin{eqnarray*}
K^k_{s} & = & - \int_{t_{k}}^{s} \biggl\{\frac{\partial\varphi
_k}{\partial t}(r,X_r^k) +
\mathcal{L}\varphi_{k}(r,X^k_{r})\\
& &\hspace*{28.8pt}{} +f\bigl(X^k_{r},\varphi_{k}(r,X^k_{r})-\alpha^*,
\sigma^{\intercal}(X^k_{r})D_x\varphi_{k}(r,X^k_{r})\bigr) \biggr\}\,dr \\
& &{} - \int_{t_{k}}^s\int_{E}( \varphi_{k}-\alpha^*
-v^*)\bigl(r,X^k_{r^-}+\gamma(X^k_{r^-},e)\bigr)\mu(dr,de) \\
& &{} + \bigl(\varphi_{k}(\theta_{k},X^k_{\theta_{k}})-\alpha^*-v(
\theta_{k},X^k_{\theta_{k}}) \bigr)\mathbf{1}_{\{s=\theta_{k}\}}.
\end{eqnarray*}
By construction and from It\^{o}'s formula on
$\varphi_k(s,X_s^k)$, we see that $(Y^k, Z^k,U^k,\break K^k)$ satisfies
(\ref{EDSRPPD}) on $[t_k,\theta_k]$. From (\ref{cond0}), it is
clear that the process $U^k$ satisfies the constraint
\[
h(U_t^k(e),e) \geq0,\qquad d\mathbf{P}\otimes dt\otimes\lambda(de),
\mbox{a.e. on } \Omega\times[t_k,\theta_k]\times E.
\]
Observe also that
%
%
\begin{equation}\label{cond2}
\varphi_{k}(\theta_{k},X^k_{\theta_{k}})-\alpha^* \geq
v(\theta_{k},X^k_{\theta_{k}}).
\end{equation}
Indeed, we have two cases:
\begin{itemize}
\item$(\theta_{k},X^k_{\theta_{k}})\in[t_{k},T]\times
B(x_{k},\delta)^c$: since $\alpha^*<\frac{\zeta}{2}$, we have by
(\ref{majhdom}),
\[
\varphi_{k}(\theta_{k},X^k_{\theta_{k}})-\alpha^* \geq
v^*(\theta_{k},X^k_{\theta_{k}}) \geq v(\theta_{k},X^k_{\theta_{k}}).
\]
%
%
\item$(\theta_{k},X^k_{\theta_{k}})\in[t_{k},T]\times
B(x_{k},\delta)\subset\mathcal{O}$: since $\alpha^* \leq\varepsilon
$, we have by (\ref{cond0})
\[
\varphi_{k}(\theta_{k},X^k_{\theta_{k}})-\alpha^* \geq
\varphi(\theta_{k},X^k_{\theta_{k}})-\varepsilon\geq g(X^k_{T})
=v(\theta_{k},X^k_{\theta_{k}}).
\]
\end{itemize}
Let us then check that $K^k$ is nondecreasing on $[t_{k},\theta_{k}]$.
First, on $[t_k,\theta_k)$, we notice that $K^k$ consists only in the Lebesgue
term $dr$, and so is nondecreasing by (\ref{cond1}). Moreover, we see
that $K^k_{\theta_{k}} \geq K^k_{\theta_{k}^-}$.
Indeed, there are two possible cases:
\begin{itemize}
\item$\theta_{k}<\tau_{k}$: then $K^k_{\theta_{k}} = K^k_{\theta
_{k}^-}+\varphi_{k}(\theta_{k},X^k_{\theta_{k}})-\alpha
^*-v(\theta_{k},X^k_{\theta_{k}})$,
and by (\ref{cond2}), we have $K^k_{\theta_{k}} \geq K^k_{\theta
_{k}^-}$.
\item$\theta_{k}=\tau_{k}$: then $K^k_{\theta_{k}} = K^k_{\theta_{k}^-}
- (\varphi_{k}(\theta_{k},X^k_{\theta_{k}})- \alpha^*-
v^*(\theta_{k},X^k_{\theta_{k}}) ) + (\varphi_{k}(\theta_{k},X^k_{\theta
_{k}})-\alpha^*-v(\theta
_{k},X^k_{\theta_{k}}) )$, and so $K^k_{\theta_{k}} \geq K^k_{\theta_{k}^-}$.
\end{itemize}
Therefore, the quadruple $(Y^k,Z^k,U^k,K^k)$ is a solution on
$[t_{k},\theta_{k}]$ to (\ref{EDSRPPD}) and (\ref{hconsPPD}), and by
Lemma \ref{lemprog}, we deduce that for all $k$,
\[
\varphi_{k}(t_{k},x_{k})-\alpha^{*} = \varphi
(t_{k},x_{k})+\sqrt{T-t_{k}}-\alpha^{*} \geq v(t_{k},x_{k}).
\]
We finally obtain a contradiction by sending $k$ to $\infty$.
\end{pf*}

\subsection{Uniqueness result}\label{sec42}

This section is devoted to a uniqueness result for the QVI (\ref
{QVIgen})--(\ref{bound-cond}).
We need to impose some additional assumptions.
\renewcommand{\theequation}{H3}
\begin{eqnarray}\label{assuH3}
\\[-19pt]
\begin{tabular}{p{325pt}}
There exists a nonnegative function
$\Lambda\in\mathcal{C}^2(\mathbb{R}^d)$ and a positive constant
$\rho$ satisfying:
\end{tabular}\hspace*{-38pt}\nonumber\\[-32pt]\nonumber
\end{eqnarray}
\begin{longlist}
\item$\mathcal{L}\Lambda+f(\cdot,\Lambda,\sigma^{\intercal}D\Lambda
)\leq\rho
\Lambda$,
\item$\inf_{e\in E}h(\mathcal{H}^e\Lambda(x)-\Lambda(x),e)>0$ for
all $x\in\mathbb{R}^d$,
\item$\Lambda(x)\geq g(x)$ for all $x\in\mathbb{R}^d$,
\item$\lim_{|x|\rightarrow\infty}\frac{\Lambda(x)}{1+|x|}=\infty$.
\end{longlist}

Assumption (\ref{assuH3}) is similar to the one made in \cite{tanyon93} or
\cite{bou06}, and essentially ensures the existence of a suitable
strict supersolution to (\ref{QVIgen}). We shall give in Section
\ref
{secsuff} some
sufficient conditions for (\ref{assuH3}). This strict supersolution allows
to control the nonlocal term in QVI (\ref{QVIgen})--(\ref{bound-cond})
via some convex
small perturbation. Thus, to deal with the dependence of $f$, $c$ on
$y,z$, we also require some convexity conditions.
\renewcommand{\theequation}{H4}
\begin{eqnarray}\label{assuH4}
\\[-19pt]\nonumber
\begin{tabular}{p{325pt}}
\mbox{ (i)} The function $f(x,\cdot,\cdot)$ is convex in $(y,z) \in\mathbb
{R}\times
\mathbb{R}^d$ for all $x \in\mathbb{R}^d$.
\end{tabular}\hspace*{-10pt}\nonumber\\[-42pt]\nonumber
\end{eqnarray}
\begin{longlist}[\hspace*{5.1pt}(iii)]
\item[(ii)] The function $h(\cdot,e)$ is concave in $u \in\mathbb{R}$ a
for all
$e \in E$.
\item[(iii)] The function $c(x,\cdot,\cdot,e)$ is convex in $(y,z) \in\mathbb
{R}\times
\mathbb{R}^d$ for all $(x,e) \in\break \mathbb{R}^d\times E$.
\item[(iv)] The function $c(x,\cdot,z,e)$ is decreasing in $y \in\mathbb
{R}$ for
all $(x,z,e) \in\mathbb{R}^d\times\mathbb{R}^d\times E$.
\end{longlist}
\begin{Theorem} \label{theouni}
Assume that (\ref{assuH3}) and (\ref{assuH4}) hold, and let $U$ (resp., $V$) be
a lsc
(resp., usc) viscosity supersolution (resp., subsolution) to (\ref
{QVIgen})--(\ref{bound-cond}) satisfying a linear growth condition
\[
\sup_{x\in\mathbb{R}^d} \frac{|U(t,x)| + |V(t,x)|}{1+|x|} <
\infty\qquad
\forall t \in[0,T].
\]
Then, $U\geq V$ on $[0, T ]\times\mathbb{R}^d$.
Consequently, under (\ref{assuH2prime}) [or (\ref{assuH1prime}) in the case: $h(u,e) = -u$],
(\ref{assuH3}), (\ref{assuH4}) and (\ref{assuHE}), the function $v$ in
(\ref{defvY}) is the unique viscosity solution to (\ref
{QVIgen})--(\ref{bound-cond}) satisfying a linear growth condition,
and $v$ is
continuous on $[0,T)\times\mathbb{R}^d$.
\end{Theorem}
\begin{pf} \textit{Comparison principle.}
As usual, we shall argue by contradiction by assuming that
%
%
\setcounter{equation}{26}
\renewcommand{\theequation}{\arabic{section}.\arabic{equation}}
\begin{equation}\label{hypabs}
\sup_{[0,T]\times\mathbb{R}^d}(V-U)> 0.
\end{equation}

1. For some $\lambda>0$ to be chosen below, let
\[
\tilde{U}(t,x) = e^{(\rho+\lambda) t}U(t,x),\qquad
\tilde{V}(t,x) = e^{(\rho+\lambda) t} V(t,x)
\]
and
\[
\tilde{\Lambda
}(t,x)=e^{(\rho+\lambda) t}\Lambda(x).
\]
A straightforward derivation shows that $\tilde{U}$ (resp., $\tilde
{V}$) is a viscosity supersolution (resp., subsolution) to
%
%
\begin{eqnarray}\label{IQVGM1}\hspace*{33pt}
\min\biggl[\rho w -\frac{\partial w}{\partial t}
-\mathcal{L}w -\tilde{f} (\cdot,w,\sigma^{\intercal}D_x w
),
\inf_{e\in E}\tilde{h} (\cdot,\tilde{\mathcal{H}}^{e} w -w
,e ) \biggr] &=& 0\nonumber\\[-8pt]\\[-8pt]
\eqntext{\mbox{on } [0,T)\times\mathbb{R}^d}\\
\label{IQVGM2}
\min\Bigl[ w(T^-,\cdot)-\tilde{g} ,\inf _{e\in E}
\tilde{h}\bigl(T,\tilde{\mathcal{H}}^{e}w(T^-,\cdot)-w(T^-,\cdot),e\bigr)\Bigr]
& = & 0 \nonumber\\[-8pt]\\[-8pt]
\eqntext{\mbox{on } \mathbb{R}^d,}
\end{eqnarray}
where
\begin{eqnarray*}
\tilde{f}(t,x,r,q)& = & e^{(\rho+\lambda) t}f \bigl(x,r e^{-(\rho
+\lambda) t},q e^{-(\rho+\lambda) t} \bigr)-\lambda r, \\
\tilde{h}(t,r,e) &= & e^{(\rho+\lambda) t}h\bigl(e^{-(\rho+\lambda)
t}r,e\bigr),\qquad \tilde{g}(x) = e^{(\rho+\lambda)T}g(x)
\end{eqnarray*}
and
\[
\tilde{\mathcal{H}}w(t,x) =
w\bigl(t,x+\gamma(x,e)\bigr)+\tilde{c} (x,w(t,x),\sigma^{\intercal}(x) D_x
w(t,x),e )
\]
with
\[
\tilde{c}(t,x,r,q,e) =
e^{(\rho+\lambda) t}c\bigl(x,e^{-(\rho+\lambda) t}r,e^{-(\rho+\lambda)
t}q,e\bigr)
\]
for all
$(t,x,r,q,e)\in[0,T]\times\mathbb{R}^d\times\mathbb{R}\times
\mathbb{R}^d\times
E$. Since $f$ is Lipschitz, we can choose $\lambda$ large enough so
that $\tilde{f}$ is nonincreasing in $r$. Denote
$\tilde{W}=(1-\mu)\tilde{U}+\mu\tilde{\Lambda}$ with $\mu>0$. By
(\ref{hypabs}) and the growth condition (\ref{assuH3})(iv) of
$\Lambda$, we have for $\mu$ small enough
%
%
\begin{equation}\label{hypabs2}
\sup_{[0,T]\times\mathbb{R}^d}(\tilde{V}-\tilde{W})=(\tilde
{V}-\tilde{W})(t_0,x_0)> 0
\end{equation}
for some $(t_0,x_0)\in[0,T]\times\mathbb{R}^d$. Moreover, from the
viscosity supersolution property (\ref{IQVGM1}) and (\ref{IQVGM2}) of
$\tilde U$, and the
conditions (\ref{assuH3})(i), (ii), (\ref{assuH4})(i), (ii), (iii), we see that
$\tilde W$ is a viscosity supersolution to
%
%
\begin{eqnarray}\label{IQV1tildeW}
\rho w -\frac{\partial w}{\partial t} -\mathcal{L}w -\tilde{f}
(\cdot,w,\sigma^{\intercal}D_x w )
& \geq& 0\qquad \mbox{on } [0,T)\times\mathbb{R}^d, \\
\label{IQV2tildeW}
\inf_{e\in E}\tilde{h} (\cdot,\tilde{\mathcal{H}}^{e} w -w
,e ) & \geq& \mu\tilde q\qquad \mbox{on }
[0,T]\times\mathbb{R}^d,
\end{eqnarray}
where $\tilde q(t,x) = e^{(\rho+\lambda)t} \inf_{e\in
E}h(\mathcal{H}^{e}\Lambda(x)-\Lambda(x),e)$ is positive on
$[0,T]\times\mathbb{R}^d$ by (\ref{assuH3})(ii).

2. Denote for all $(t,x,y)\in[0,T]\times\mathbb
{R}^d\times\mathbb{R}^d$ and $n\geq1$
\[
\Theta_n(t,x,y) = \tilde{V}(t,x)-\tilde{W}(t,y) -\varphi_n(t,x,y)
\]
with
\[
\varphi_n(t,x,y) = n|x-y|^2+|x-x_0|^4+|t-t_0|^2.
\]
By the growth assumption on $U$ and $V$ and (\ref{assuH3})(iii), for all
$n$, there exists
$(t_{n},x_{n},y_{n})\in[0,T]\times\mathbb{R}^d\times\mathbb{R}^d$
attaining the maximum of $\Theta_n$ on $[0,T]\times\mathbb
{R}^d\times\mathbb{R}
^d$. By
standard arguments, we have
%
%
\begin{eqnarray} \label{limtn}
(t_n,x_n,y_n) & \rightarrow& (t_0,x_0,x_0), \\
\label{limvarphin}
n|x_n-y_n|^2 & \rightarrow& 0, \\
\label{limvwn}
\tilde V(t_n,x_n)-\tilde W(t_n,y_n) &\rightarrow& \tilde
V(t_0,x_0)-\tilde W(t_0,x_0).
\end{eqnarray}

3. We now show that for $n$ large enough
%
%
\begin{equation}\label{sursoltildeV}
\inf_{e\in
E}\tilde{h}\bigl(t_n,\tilde{\mathcal{H}}^{e}[t_n,x_n,D_x\varphi
_{n}(t_n,x_n,y_n),\tilde{V}]-\tilde{V}(t_{n},x_{n}),e\bigr) > 0.
\end{equation}
On the contrary, up to a subsequence, we would have for all $n$,
\[
\inf_{e\in E}
\tilde{h}\bigl(t_n,\tilde{\mathcal{H}}^{e}[t_n,x_n,D_x\varphi
_{n}(t_n,x_n,y_n),\tilde{V}]-\tilde{V}(t_{n},x_{n}),e\bigr) \leq0
\]
and so by uppersemicontinuity of $\tilde V$, compactness of $E$, there
would exist a sequence $(e_n)$ in $E$ such that
\[
\tilde{h}\bigl(t_n,\tilde{\mathcal{H}}^{e_n}[t_n,x_n,D_x\varphi
_{n}(t_n,x_n,y_n),\tilde{V}]-\tilde{V}(t_{n},x_{n}),e_n\bigr) \leq0.
\]
Moreover, by the viscosity supersolution property of $\tilde W$ to
(\ref{IQV2tildeW}), we have
\[
\tilde{h}\bigl(t_n,\tilde{\mathcal{H}}^{e_n}[t_n,y_n,-D_y\varphi
_{n}(t_n,x_n,y_n),\tilde{W}] -
\tilde{W}(t_n,y_n),e_n\bigr) \geq\mu\tilde{q}(t_n,y_n).
\]
From the nonincreasing and the Lipschitz property of $h(\cdot,e)$, we
deduce from the two previous ine\-qualities that there exists a
positive constant $\eta$ such that
\begin{eqnarray*}
&& \tilde{\mathcal{H}}^{e_{n}}[t_n,y_n,-D_y\varphi
_{n}(t_n,x_n,y_n),\tilde{W}] - \tilde{W}(t_n,y_n)+\eta\tilde
{q}(t_n,y_n) \\
&&\qquad \leq \tilde{\mathcal{H}}^{e_{n}}[t_n,x_n,D_x\varphi
_{n}(t_n,x_n,y_n),\tilde{V}]-\tilde{V}(t_n,x_n),
\end{eqnarray*}
which is rewritten as
%
%
\begin{eqnarray} \label{intervwn}
&& \tilde{V}(t_n,x_n)-\tilde{W}(t_n,y_n)+\eta\tilde{q}(t_n,y_n)
\nonumber\\[-8pt]\\[-8pt]
&&\qquad\leq \tilde{V}\bigl(t_{n},x_n+\gamma(x_n,e_n)\bigr) -\tilde
{W}\bigl(t_n,y_n+\gamma(y_n,e_n)\bigr) +\Delta C_n,\nonumber
\end{eqnarray}
%
where
\begin{eqnarray*}
\Delta C_n &= & \tilde{c} (t_{n},x_n,\tilde{V}(t_n,x_n),\sigma
^{\intercal}(x_{n}) D_x\varphi_{n}(t_n,x_n,y_n),e_n )\\
& &{} -
\tilde{c} (t_{n},y_n,\tilde{W}(t_n,y_n),-\sigma^{\intercal}(y_{n})
D_y\varphi_{n}(t_n,x_n,y_n) ).
\end{eqnarray*}
Now, we write $\Delta C_n = \Delta C_n^1 + \Delta C_n^2 + \Delta
C_n^3$, with
\begin{eqnarray*}
\Delta C_n^1 &=& \tilde{c} (t_n,x_n,\tilde{V}(t_n,x_n),\sigma
^{\intercal}(x_{n}) D_x\varphi_n(t_n,x_n,y_n),e_n ) \\
& &{} - \tilde{c} (t_n,x_n,\tilde{W}(t_n,y_n),\sigma
^{\intercal}
(x_{n}) D_x\varphi_n(t_n,x_n,y_n),e_n ), \\
\Delta C_n^2 &=& \tilde{c} (t_n,x_n,\tilde{W}(t_n,y_n),\sigma
^{\intercal}(x_{n}) D_x\varphi_n(t_n,x_n,y_n),e_n ) \\
& &{} - \tilde{c} (t_n,x_n,\tilde{W}(t_n,y_n),-\sigma
^{\intercal}
(y_{n}) D_y\varphi_{n}(t_n,x_n,y_n),e_n ), \\
\Delta C_n^3 &=&
\tilde{c} (t_n,x_n,\tilde{W}(t_n,y_n),-\sigma^{\intercal}(y_{n})
D_y\varphi_{n}(t_n,x_n,y_n),e_n ) \\
& &{} -
\tilde{c} (t_n,y_n,\tilde{W}(t_n,y_n),-\sigma^{\intercal}(y_{n})
D_y\varphi_{n}(t_n,x_n,y_n),e_n ).
\end{eqnarray*}
We have $\tilde
V(t_n,x_n)-\tilde W(t_n,y_n) \rightarrow(\tilde V-\tilde
W)(t_0,x_0) > 0$ by (\ref{hypabs2}) and (\ref{limvwn}). Hence,
for $n$ large enough, $\tilde V(t_n,x_n) \geq\tilde
W(t_n,y_n)$, and so from the\break nonincreasing condition
(\ref{assuH4})(iv) of $c$, we have $\Delta C_n^1 \leq0$. Since\break
$\sigma^{\intercal}(x_{n})D_{x}\varphi_{n}(t_{n},x_{n},y_{n})+
\sigma^{\intercal}(y_{n})D_{y}\varphi_{n}(t_{n},x_{n},y_{n})
\rightarrow0$ by the Lipschitz condition on $\sigma$ and
(\ref{limvarphin}), we deduce with the Lipschitz condition on $c$
that $\limsup_{n\rightarrow\infty} \Delta C_n^2 \leq0$. By
(\ref{limtn}) and continuity of $c$, we have
$\lim_{n\rightarrow\infty} \Delta C_n^3 = 0$. Therefore, we
obtain
\[
\limsup_{n\rightarrow\infty} \Delta C_n \leq0.
\]
Up to a subsequence, we may assume that $(e_n)$ converges to
$e_0$ in $E$. Hence, by sending $n$ to infinity into
(\ref{intervwn}), it follows with (\ref{limvwn}) and the upper
(resp., lower)-semicontinuity of $\tilde V$ (resp., $\tilde W$)
that
\begin{eqnarray*}
(\tilde V-\tilde W)\bigl(t_0,x_0+\gamma(x_0,e_0),x_0+\gamma(x_0,e_0)\bigr) &
\geq& (\tilde V-\tilde W)(t_0,x_0) + \eta\tilde{q}(t_{0},x_0) \\
& > & (\tilde V-\tilde W)(t_0,x_0),
\end{eqnarray*}
a contradiction with (\ref{hypabs2}).

4. Let us check that, up to a subsequence,
$t_n<T$ for all $n$. On the contrary, $t_n = t_0 = T$ for $n$
large enough, and
from (\ref{sursoltildeV}), and the viscosity subsolution property of
$\tilde V$ to (\ref{IQVGM2}), we would get
\[
\tilde{V}(T,x_{n}) \leq\tilde{g}(x_{n}).
\]
On the other hand, by the viscosity supersolution property of $\tilde
U$ to (\ref{IQVGM2}) and (\ref{assuH3})(iii), we have $\tilde W(T,y_n) \geq
\tilde g(y_n)$,
and so
\[
\tilde V(T,x_n) -\tilde W(T,y_n) \leq\tilde g(x_n)- \tilde g(y_n).
\]
By sending $n$ to infinity, and from continuity of $\tilde g$, this
would imply $(\tilde V-\tilde W)(t_0,x_0) \leq0$, a contradiction
with (\ref{hypabs2}).

5. We may then apply Ishii's lemma (see Theorem
6.1 in \cite{fleson06}) to $(t_n,x_n,y_n) \in[0,T)\times\mathbb{R}
^d\times\mathbb{R}^d$
that attains the maximum of $\Theta_n$, for all $n \geq1$:
there exist $(p_{\tilde{V}}^n,q_{\tilde{V}}^n,M_n)\in\bar
{J}^{2,+}\tilde{V}(t_n,x_n)$ and
$(p_{\tilde{W}}^n,q_{\tilde{W}}^n,N_n)\in\bar{J}^{2,-}\tilde
{W}(t_n,y_n)$ such that
\begin{eqnarray*}
p_{\tilde{V}}^n-p_{\tilde{W}}^n &=& \partial_t\varphi
_n(t_n,x_n,y_n) = 2(t_n-t_0),\\
q_{\tilde{V}}^n &=& D_x\varphi_n(t_n,x_n,y_n), \qquad q_{\tilde
{W}}^n = -D_y\varphi_n(t_n,x_n,y_n)
\end{eqnarray*}
and
%
%
\begin{equation}\label{condishii}
\pmatrix{M_{n} & 0 \cr0 & -N_{n}}
\leq A_n+ \frac{1}{2n} A_n^2,
\end{equation}
where
$A_{n}=D^2_{(x,y)}\varphi_n(t_n,x_n,y_n)$. From the viscosity
supersolution property of $\tilde W$ to (\ref{IQV1tildeW}), we have
\begin{eqnarray*}
&&\rho\tilde{W}(t_n,y_n)-p_{\tilde{W}}^n+\langle b(y_n), D_y\varphi
(t_n,x_n,y_n)\rangle-
\tfrac{1}{2}\operatorname{tr} (\sigma(y_n)\sigma^{\intercal}(y_n)
N_n )
\\
&&\qquad{}
-\tilde{f} (t_n, y_n,\tilde{W}(t_n,y_n),-\sigma^{\intercal}(y_{n})
D_y\varphi(t_n,x_n,y_n) ) \geq 0.
\end{eqnarray*}
On the other hand,
from (\ref{sursoltildeV}) and the viscosity subsolution property
of $\tilde V$ to (\ref{IQVGM1}), we have
\begin{eqnarray*}
&&\rho
\tilde{V}(t_n,x_n)-p_{\tilde{V}}^n - \langle b(x_n),
D_x\varphi(t_n,x_n,y_n)\rangle
-\tfrac{1}{2}\operatorname{tr} (\sigma(x_n)\sigma^{\intercal}(x_n)
M_n )
\\
&&\qquad{} - \tilde{f} (t_n,
x_n,\tilde{V}(t_n,x_n),\sigma^{\intercal}(x_{n})D_x\varphi
(t_n,x_n,y_n) )
\leq 0.
\end{eqnarray*}
By subtracting the two previous inequalities, we obtain
%
%
\begin{eqnarray}\label{interfin}
&&\rho\bigl(\tilde V(t_n,x_n)-\tilde W(t_n,y_n)\bigr)\nonumber\\
&&\qquad \leq p_{\tilde{V}}^n -
p_{\tilde{W}}^n + \Delta F_n \nonumber\\[-8pt]\\[-8pt]
&&\qquad\quad{} + \langle b(x_n),D_x\varphi_n(t_n,x_n,y_n)\rangle+ \langle b(y_n),
D_y\varphi_n(t_n,x_n,y_n) \rangle\nonumber\\
&&\qquad\quad{} +
\tfrac{1}{2}\operatorname{tr} \bigl(\sigma(x_n)\sigma^{\intercal
}(x_n)M_n-\sigma
(y_n)\sigma^{\intercal}(y_n)N_n \bigr),\nonumber
\end{eqnarray}
where
\begin{eqnarray*}
\Delta F_n &=& \tilde{f}(t_n,x_n,\tilde{V}(t_n,x_n),\sigma
^{\intercal}
(x_{n})D_x\varphi_n(t_n,x_n,y_n)) \\
& &{} -
\tilde{f}(t_n,y_n,\tilde{W}(t_n,y_n),-\sigma^{\intercal
}(y_{n})D_y\varphi
_n(t_n,x_n,y_n)).
\end{eqnarray*}
From (\ref{limtn}), we have
$p_{\tilde{V}}^n-p_{\tilde{W}}^n \rightarrow0$ as $n$ goes
to infinity. From the Lipschitz property of $b$, and
(\ref{limvarphin}), we have
\[
\lim_{n\rightarrow\infty} \bigl(\langle b(x_n), D_x\varphi
_n(t_n,x_n,y_n)\rangle+
\langle b(y_n), D_y\varphi_n(t_n,x_n,y_n)
\bigr) = 0.
\]
As usual, from (\ref{condishii}), (\ref{limtn}),
(\ref{limvarphin}) and the Lipschitz property of $\sigma$, we have
\[
\limsup_{n\rightarrow\infty} \operatorname{tr}
\bigl(\sigma
(x_n)\sigma^{\intercal}(x_n)M_n-\sigma(y_n)
\sigma^{\intercal}(y_n)N_n \bigr)
\leq0.
\]
Moreover, by the same arguments as for $\tilde c$,
using the nonincreasing property of $\tilde f$ in its third
variable, and the Lipschitz property of $\tilde f$, we have
\[
\limsup_{n\rightarrow\infty} \Delta F_n \leq0.
\]
Therefore, by
sending $n\rightarrow\infty$ into (\ref{interfin}), we conclude with
(\ref{limvwn}) that $\rho(\tilde V-\tilde W)(t_0,x_0) \leq0$, a
contradiction with (\ref{hypabs2}).

\textit{Uniqueness for $v$}. The uniqueness result is
then a direct consequence of the comparison principle, and the
continuity of $v$ on $[0,T)\times\mathbb{R}^d$ follows from the fact
that in
this case $v_* = v^*$.
\end{pf}
\begin{Remark} \label{remconvcont}
As a byproduct of the comparison principle in Theorem \ref
{theouni}, we get the continuity of the value function $v$ on
$[0,T)\times\mathbb{R}^d$. Since the jump-diffusion process $X$ is quasi-left
continuous, then so is the minimal solution
$Y_t = v(t,X_t)$ to the BSDE with constrained jumps, and the
penalized approximation $Y_t^n = v_n(t,X_t)$. This implies that the
predictable
projections $^{p}Y$ and $^{p}Y^n$, respectively, of $Y$ and $Y^n$, are
equal to $^{p}Y_t = Y_{t^-}$ and $^{p}Y_t^n = Y_{t^-}^n$. Therefore,
$Y_{t^-}=\lim_{n\to\infty} Y_{t^-}^n$. From the weak version of
Dini's theorem (see \cite{delmey80}, page 202) this yields the uniform
convergence of
$Y^n$ on $[0,T]$, that is, $\lim_{n\rightarrow\infty} \sup_{t\in[0,T]}
|Y_t^n-Y_t| = 0$, and so by the dominated convergence theorem, the
convergence of $Y^n$ to $Y$ in $\bolds{\mathcal{S}}^{\mathbf{2}}$:
%
%
\setcounter{equation}{39}
\renewcommand{\theequation}{\arabic{section}.\arabic{equation}}
\begin{equation}\label{limYnfort}
{\lim_{n\rightarrow\infty}} \|Y^n -Y\|_{{\bolds{\mathcal{S}}^{\mathbf
{2}}}} = 0.
\end{equation}
%
Then, by applying It\^{o}'s formula to $t\mapsto\mathbf
{E}|Y_{t}-Y_{t}^n|$ a
in the proof of Theorem \ref{thmmain1}, we get from the convergence of
$Y^n$ to $Y$ in $\bolds{\mathcal{S}}^{\mathbf{2}}$ that $(Z^n,V^n)$
converges to
$(Z,V)$ in
$\mathbf{L}^{\mathbf{2}}(\mathbf{W})\times\mathbf{L}^{\mathbf{2}}(\tilde
\mu)$ and that $K$ is continuous.
\end{Remark}

\section[Some sufficient conditions for (H2') and (H3)]{Some sufficient
conditions for (\protect\ref{assuH2prime}) and (\protect\ref{assuH3})}
\label{secsuff}

In this section, we provide various explicit conditions on the
coefficients model, which ensure that the general assumptions
(\ref{assuH2prime}) and (\ref{assuH3}) hold true.

\subsection{Existence of the solution to BSDE with jump
constraint}\label{sec51}

We first consider a case where we have upper bounds for the
coefficients and $h(u,e)=-u$.
\begin{Proposition} Suppose that $h(u,e) = -u$, and assume that
there exist real constants $C_{1},C_{2}$ and $\eta\in\mathbb{R}^{d}$
such that
%
%
\setcounter{equation}{0}
\renewcommand{\theequation}{\arabic{section}.\arabic{equation}}
\begin{eqnarray}\quad
\label{hypcons}
g(x) &\leq& C_{1}+\langle\eta,
x\rangle,\nonumber\\[-8pt]\\[-8pt]
c(x,y,z,e)+\langle\eta,\gamma(x,e)\rangle
&\leq& 0\quad\mbox{and}\quad
f(x,y,z)+\langle\eta, b(x)\rangle \leq C_2 \nonumber
\end{eqnarray}
for all $(x,y,z,e)\in\mathbb{R}^{d}\times\mathbb{R}\times\mathbb
{R}^{d}\times E$. Then (\ref{assuH2prime}) holds true.
\end{Proposition}
\begin{pf} Let us define a quadruple $(\tilde Y,\tilde Z,\tilde
K,\tilde U)$ by: $\tilde Y_t=C_1+C_2(T-t)+\langle\eta, X_t\rangle
$ for
$t < T$, $\tilde Y_T = g(X_T)$, $\tilde
Z_t=\sigma(X_{t^-}).\eta$, $\tilde U_t(e)=0$ and
\begin{eqnarray*}
\tilde K_t &= & \int_0^t \{C_2-\eta\cdot b(X_s)-f(X_s,\tilde
Y_s,\tilde Z_s) \}\,ds \\
& &{} - \int_0^t\int_E \{c(X_{s-},\tilde Y_{s-},\tilde
Z_s,e)+\langle\eta,\gamma(X_{s-},e)\rangle
\}\mu(ds,de),\qquad t < T, \\
\tilde K_T & = & \tilde K_{T-}+C_1+\langle\eta, X_T\rangle-g(X_T).
\end{eqnarray*}
From (\ref{hypcons}), the process $\tilde K$ is clearly
nondecreasing. Moreover, from the dynamics of $X$, and by
construction, we see that the quadruple $(\tilde Y,\tilde Z,\tilde
K,\tilde U)$ satisfies (\ref{BSDEgen})--(\ref{negjump}) and the
function $\tilde v(t,x) = C_1+C_2(T-t)+\eta.x$ clearly satisfies
a linear growth condition.
\end{pf}

We next give an example inspired by \cite{bou06} where the jumps
of $X$ vanish as $X$ goes out of a ball centered in zero in the case
of impulse control.
\begin{Proposition} Suppose that $h(u,e) = -u$, $f,c$ does not
depend on $y,z$, and assume that $c\leq0$, $\gamma=0$ on
$\{x\in\mathbb{R}^{d}\dvtx|x|\geq C_1 \}\times E$ for some $C_1>0$.
Then, (\ref{assuH2prime}) holds true.
\end{Proposition}
\begin{pf}
We consider the function $v$
\[
v(t,x) = \sup_{\nu\in\mathcal{V}}\mathbf{E}^{\nu} \biggl[
g(X_{T}^{t,x})+\int
_{t}^Tf(X_{s}^{t,x})\,ds+\int_{t}^T\int_{E}c(X_{s^-}^{t,x},e)\mu
(ds,de) \biggr].
\]
Since $c \leq0$, and the choice of $\nu= 1$ corresponds to
the probability measure $\mathbf{P}^1 = \mathbf{P}$, we see that
$\hat v \leq v \leq\bar v$ where
\begin{eqnarray*}
\hat v(t,x) &=& \mathbf{E} \biggl[ g(X_{T}^{t,x})+\int
_{t}^Tf(X_{s}^{t,x})\,ds+\int_{t}^T\int_{E}c(X_{s^-}^{t,x},e)\mu
(ds,de) \biggr], \\
\bar v(t,x) &=& \sup_{\nu\in\mathcal{V}}\mathbf{E}^{\nu} \biggl[
g(X_{T}^{t,x})+\int_{t}^Tf(X_{s}^{t,x})\,ds \biggr].
\end{eqnarray*}
The function $\hat v$ clearly satisfies a linear growth condition by
the linear growth conditions on $g,f,c$ and the standard estimate for
$X$. Moreover,
under the assumptions on the jump coefficient $\gamma$, it is shown in
\cite{bou06} that $\bar v$ satisfies a linear growth condition.
Therefore, $\hat v$ also satisfies a linear growth condition.

Let us now define the process $Y_t = v(t,X_t)$, which is then
equal to
\[
Y_t = \mathop{\operatorname{ess}\sup}_{\nu\in\mathcal{V}}\mathbf
{E}^{\nu} \biggl[
g(X_{T})+\int_{t}^Tf(X_{s})\,ds+\int_{t}^T\int_{E}c(X_{s^-},e)\mu(ds,de)
\Big| \mathcal{F}_t \biggr],
\]
and lies in $\bolds{\mathcal{S}}^{\mathbf{2}}$ from the linear
growth condition, and the estimate (\ref{Xest}) for $X$. From
Theorem \ref{theonegjump}, we then know that there exists $(Z,U,K) \in
\mathbf{L}^{\mathbf{2}}(\mathbf{W})\times\mathbf{L}^{\mathbf{2}}(\tilde
\mu)\times\mathbf{A}^{\mathbf{2}}$ such
that $(Y,Z,U,K)$ is the minimal solution to
(\ref{BSDEgen})--(\ref{negjump}), and so (\ref{assuH2prime}) is satisfied.
\end{pf}

We finally consider a case for general constraint function $h$.
\begin{Proposition}
Assume that there exists a Lipschitz function $w\in
\mathcal{C}^{2}(\mathbb{R}^{d})$ satisfying a linear growth
condition, supersolution to (\ref{bound-cond}), and such that
\[
\langle b, Dw\rangle+ \tfrac{1}{2}\operatorname{tr}(\sigma
\sigma^{\intercal}
D^{2}w)+ f(\cdot,w,\sigma^{\intercal}Dw) \leq C\qquad
\mbox{on }
\mathbb{R}^d
\]
for some constant $C$. Then (\ref{assuH2prime}) holds true.
\end{Proposition}
\begin{pf}
Let us define a quadruple $(\tilde Y,\tilde Z,\tilde
U,\tilde K)$ by
\[
\tilde Y_{t} = w(X_{t})+C(T-t),
\qquad t < T,\qquad \tilde Y_{T}=g(X_{T}),
\]
$\tilde
Z_{t}=\sigma^{\intercal}(X_{t^-}) Dw(X_{t^-}),\tilde
U_{t}(e)=w(X_{t^-}+\gamma(X_{t^-},e))+c(X_{t^-},\tilde
Y_{t^-},\tilde Z_t,e)-w(X_{t^-})$, and
\begin{eqnarray*}
\tilde K_t & = &
\int_0^t \biggl[C- \langle b(X_s), Dw(X_{s})\rangle
\\
& &\hspace*{16.85pt}{} -\frac{1}{2}\operatorname{tr}\{\sigma(X_{s})\sigma^{\intercal}(X_{s})
D^{2}w(X_{s})\}-f(X_s,\tilde Y_s,\tilde Z_s) \biggr]\,ds,\qquad
t < T, \\
\tilde K_T &= & \tilde K_{T^-} + w(X_{T})-g(X_{T}).
\end{eqnarray*}
From the
conditions on $w$, we see that $(\tilde Y,\tilde Z,\tilde K,\tilde
U)$ lies in $\bolds{\mathcal{S}}^{\mathbf{2}}\times\mathbf{L}^{\mathbf
{2}}(\mathbf{W})\times{\bf
L^2(\tilde
\mu)}\times\mathbf{A}^{\mathbf{2}}$. Moreover, by It\^{o}'s formula to
$w(X_t)$ and
the supersolution property of $w$ to~(\ref{bound-cond}), we
conclude that $(\tilde Y,\tilde Z,\tilde K,\tilde U)$ is solution to
(\ref{BSDEgen}) and (\ref{hcons}), and $\tilde v(t,x) = w(t,x)+C(T-t)$
satisfies a linear growth condition.
\end{pf}

\subsection[The strict supersolution condition (H3)]{The strict supersolution
condition (\protect\ref{assuH3})}\label{sec52}

We give a sufficient condition for (\ref{assuH3}) in the usual case where
$f$ and $c$ do not depend neither on $y$ nor on $z$.
\begin{Proposition}
Consider the case where $h$ is given by
\[
h(u,e)=-u.
\]
Assume that there exists a constant $\alpha>0$ such that
\begin{eqnarray*}
-\alpha& < & |x+\gamma(x,e)|^2-|x|^2\qquad \forall(x,e)\in\mathbb{R}
^d\times E,\\
\beta :\!&= & \inf_{(x,e)\in\mathbb{R}^d\times E}\frac
{-c(x,e)}{|x+\gamma
(x,e)|^2-|x|^2+\alpha}>0.
\end{eqnarray*}
Then assumption (\ref{assuH3}) holds true.
\end{Proposition}
\begin{pf}
We set $\Lambda(x):=\beta|x|^2+\zeta$ with $\zeta$ large enough so
that $\Lambda\geq g$, that is, (\ref{assuH3})(iii) is satisfied.
A straightforward computation shows that
\[
\inf_{e\in E}h\bigl(\mathcal{H}^{e}\Lambda(x)-\Lambda(x),e\bigr)\geq\alpha
\beta>0
\]
and hence (\ref{assuH3})(ii) is satisfied. Clearly, (\ref{assuH3})(iv)
holds as well. Finally, it follows from the linear growth assumption
on $b$ and $\sigma$ that (\ref{assuH3})(i) holds for a sufficiently
large parameter $\rho$.
\end{pf}

\section*{Acknowledgments}
We would like to thank both referees for useful comments.

%

%
\printaddresses

\end{document}